\documentclass{article}
\usepackage{amsmath,amsthm}
\usepackage{graphicx,epstopdf,epsfig,multirow,epic,bm}
\usepackage{subfigure}
\usepackage{amsfonts}
\usepackage{amssymb}
\usepackage{lineno,hyperref}

\normalsize \rm
\begin{document}

\begin{center}
{\Large  \textbf {Random walks on stochastic uniform growth trees: Analytical formula for mean first-passage time}}\\[12pt]
{\large \quad Fei Ma$^{a,}$\footnote{~The author's E-mail: mafei123987@163.com. } \; and \; Ping Wang$^{b,c,d}$}\\[6pt]
{\footnotesize $^{a}$ School of Electronics Engineering and Computer Science, Peking University, Beijing 100871, China\\
$^{b}$ National Engineering Research Center for Software Engineering, Peking University, Beijing, China\\
$^{c}$ School of Software and Microelectronics, Peking University, Beijing  102600, China\\
$^{d}$ Key Laboratory of High Confidence Software Technologies (PKU), Ministry of Education, Beijing, China
}\\[12pt]
\end{center}

\begin{quote}
\textbf{Abstract:}  As known, the commonly-utilized ways to determine mean first-passage time $\overline{\mathcal{F}}$ for random walk on networks are mainly based on Laplacian spectra. However, methods of this type can become prohibitively complicated and even fail to work when the Laplacian matrix of network under consideration is difficult to describe in the first place. In this paper, we propose an effective approach to determining quantity $\overline{\mathcal{F}}$ on some widely-studied tree networks. To this end, we first build up a general formula between Wiener index $\mathcal{W}$ and $\overline{\mathcal{F}}$ on a tree. This enables us to convert issues to answer into calculation of $\mathcal{W}$ on networks in question. As opposed to most of previous work focusing on deterministic growth trees, our goal is to consider stochastic case. Towards this end, we establish a principled framework where randomness is introduced into the process of growing trees. As an immediate consequence, the previously published results upon deterministic cases are thoroughly covered by formulas established in this paper. Additionally, it is also straightforward to obtain Kirchhoff index on our tree networks using the proposed approach. Most importantly, our approach is more manageable than many other methods including spectral technique in situations considered herein. \\


\end{quote}

\section{Introduction}

Random walk on networks is a hot-topic in the realm of complex network study \cite{Noh-2004}-\cite{Lam-2012} because of a large number of applications, both theoretical and practical, such as, the study of transport-limited reactions \cite{O-B-2011}, target search \cite{Michael-2006} as well as disease spreading on relationship networks among individuals \cite{Jia-2018}, to name just a few. Thus, it has attracted considerable attention from various fields including applied mathematics, theoretical computer science, statistic physics in the past years \cite{Carletti-2020}-\cite{Berenbrink-2010}. In the language of mathematics, random walk on networks $\mathcal{G}(\mathcal{V},\mathcal{E})$ (defined in detail later) describes a simple dynamic process where a walker starting out from its current position $u\in\mathcal{V}$ will at random hop to any vertex $v\in\mathcal{V}$ in the next step with probability $P_{u\rightarrow v}$ that is defined as follows

$$P_{u\rightarrow v}=\left\{\begin{split}&1/k_{_{u}} \qquad \text{if $u$ is adjacent to $v$},\\
&0 \qquad \qquad\text{otherwise.}
\end{split}\right.$$
in which $k_{_{u}}$ is the number of edges incident with vertex $u$. In other words,

$$f(v)=\sum_{u,u\sim v}\frac{1}{k_{u}}f(u)$$
for any distribution $f:\mathcal{V}\rightarrow \mathbb{R}$ with $\sum_{v}f(v)=1$ where symbol $u\sim v$ indicates that vertex $u$ is adjacent to $v$ \cite{Chung-1998}. Along this line of research, we keep on studying random walks on some networks of significant interest by considering many relevant topological parameters in order to better understand how the underlying structure affects dynamic behaviors of such type. As a tip, the terms graph and network are used indistinctly throughout this paper.

In fact, there is a long history of investigating random walks on an arbitrary graph \cite{Shuji-2015}-\cite{Ibe-2013}. As known, the most fundamental in studying random walks is to estimate some structural parameters, for instance, mean first-passage time $\overline{\mathcal{F}}$ (defined in detail later). In theory, the closed-form solution to mean first-passage time $\overline{\mathcal{F}}$ on graph $\mathcal{G}(\mathcal{V},\mathcal{E})$ may be obtained using

$$\overline{\mathcal{F}}=\frac{2|\mathcal{E}|}{|\mathcal{V}|-1}\sum_{i=2}^{|\mathcal{V}|}\frac{1}{\lambda_{i}}$$
where $|\mathcal{E}|$ and $|\mathcal{V}|$ represents size and order of graph $\mathcal{G}(\mathcal{V},\mathcal{E})$ in the jargon of graph theory \cite{Bondy-2008}, respectively, and $\lambda_{i}$ is eigenvalue of the corresponding Laplacian matrix $\mathbf{L}_{\mathcal{G}}$. On the other hand, this typical method will become intractable to derive what we want when it is not easy to determine Laplacian spectra of network in question. Perhaps, one reason for this is that it is difficult to first establish the corresponding Laplacian matrix. These such example networks can be ubiquitous in research community, for instance, the famous Vicsek fractal \cite{Vicsek-1983}. This triggers the relevant research and inspires scholars to develop effective approaches, which are suitable for different types of networks under considerable, in order to get around the dilemma above. An instructive example is the so-called self-similarity based method for heterogeneous networks \cite{Sheng-2019}. Therefore, in this paper, we will build up feasible techniques for the purpose of determining the analytic solution to mean first-passage time on some networks, in particular, tree networks. In general, these tree models that will be discussed are significantly difficult to analyze according to the typical methods when considering random walks in more general cases as will be shown later.

Tree, as the simplest connected network, has also been well discussed in the area of random walk study \cite{Graham-1990}-\cite{Ma-2020-1}, and some structural parameters have been reported. For instance, the mean first-passage time $\overline{\mathcal{F}}_{\mathcal{T}}$ on an arbitrary tree $\mathcal{T}$ complies with the next expression

$$2|\mathcal{V}_{\mathcal{T}}|\leq\overline{\mathcal{F}}_{\mathcal{T}}\leq\frac{1}{3}|\mathcal{V}_{\mathcal{T}}|^{2}.$$
Note that we have neglected other negligible terms in the above inequality in the large graph size limit. Nonetheless, the problem of how to analytically calculate the concrete formula for mean first-passage time on trees with intriguing structural features is of considerable interest both theoretically and experimentally. This is because there are a great number of real-world applications associated with some specific trees \cite{Borah-2014}-\cite{Bartolo-2016}. As a topological measure, the mean first-passage time can be chosen to quantify the structural properties including network robustness on those tree models. Particularly, the well-known Vicsek fractal is often used to describe the underlying structure of some dendrimers and regular hyperbranched polymers \cite{Borah-2014}. The deterministic uniform growth tree is frequently adopted to model epidemic spreading in population \cite{Moon-1974}. The famous $T$-graph \cite{Redner-2001} and other variants including Peano basin fractal \cite{Bartolo-2016} have been popularly utilized in physics and geosciences. Motivated by this, we aim at studying random walks on stochastic uniform growth trees, which are a class of more general tree models, and then derive the corresponding solution to mean first-passage time. To this end, we first propose a principled framework for generating anticipated tree models. As an immediate result, those tree models mentioned above will be grouped into the proposed framework. In other wards, the goal of this study is to try to discuss those trees from a more comprehensive and systematical viewpoint.

In view of applications of those tree models into wide ranges of science \cite{Borah-2014}-\cite{Bartolo-2016}, the exact formulae for mean first-passage time on them have been obtained using some methods which are mainly based on spectral theory \cite{Kemeny-1976,Biggs-1974}. It is worth noticing that those models previously discussed all share deterministic structure. Even so, it is not easy to calculate mean first-passage time by using the typical methods as mentioned above. Meanwhile, in most cases, almost all seeds used to create those tree models are some specific and simple trees, for instance, an edge or a star. If an arbitrary tree is selected to serve as the seed, it is clear to understand that the pre-existing methods might not be suitable for calculating analytic solution to mean first-passage time. Additionally, when introducing randomness into development of tree models of this kind, the typical techniques will be prohibitively difficult and even fail to work. To address this issue, we develop some more convenient combinatorial approaches, and finally obtain what we are seeking for.

Our contribution is shown in the following form.

(1) A principled framework for creating stochastic uniform growth tree is proposed. More concretely, it consists of three ingredients: Vertex-based uniform growth mechanism, Edge-based uniform growth mechanism, along with mixture uniform growth mechanism. As will be clear to the eye, some previous models are contained in this framework completely. Meanwhile, based on this framework, we uncover some close relationships between those pre-existing models, which are previously unknown.

(2) Using a combinatorial method, we derive a formula connecting Wiener index $\mathcal{W}$ to mean first-passage time $\overline{\mathcal{F}}$ in a tree. With this formula, the analytic solutions to mean first-passage time for random walks on all stochastic uniform growth trees built in this work are obtained precisely. As a consequence, many early published results in the literature are certainly covered by the formulae derived by us. More importantly, no complicated computations are involved in the process of determining solutions.

(3) With the concept of network criticality, we find that the underlying topological structure of stochastic uniform growth tree plays a key role in studying random walks on it. At the same time, the scaling relations between quantity $\overline{\mathcal{F}}$ and vertex number on all stochastic uniform growth trees proposed herein are also analytically obtained. They can fairly be thought of as more general consequences in comparison with the previously reported results.

The rest of this paper is organized as follows. We in Section 2 introduce some conventional notations in the jargon of graph theory. In Section 3, we propose some widely-studied graphic operations to establish a principled framework for generating potential candidate models as the objectives of this paper. Some example networks are listed for the purpose of better understanding the proposed framework. Next, the main results are shown in Section 4. More specifically, a formal connection of Wiener index $\mathcal{W}$ to mean first-passage time $\overline{\mathcal{F}}$ in a tree is built in a combinatorial manner. This enables us to derive the analytic solutions to mean first-passage time for random walks on all the above-mentioned models in a more convenient way. Accordingly, some detailed discussions are provided in order to understand how the underlying structure affects dynamic behaviors of this kind. Finally, we draw conclusion in the last section.

\section{Definitions and notations}

It is a convention in graph theory to denote a graph by $\mathcal{G}(\mathcal{V},\mathcal{E})$ that contains a set $\mathcal{V}$ vertices and a set $\mathcal{E}$ edges running between vertices. Accordingly, symbols $|\mathcal{V}|$ and $|\mathcal{E}|$ represent vertex number and edge number, respectively. $\mathcal{N}_{u}$ is used to indicate the neighboring set of vertex $u$. At the same time, let notation $[a,b]$ be a collection of integers which are certainly both no greater than $b$ and no less than $a$.

\subsection{Matrix representation of graph}

Roughly speaking, it is convenient to interpret a graph $\mathcal{G}(\mathcal{V},\mathcal{E})$ based on its adjacency matrix $\mathbf{A}=(a_{uv})$, which is defined in the following form

$$a_{uv}=\left\{\begin{aligned}&1, \quad\text{vertex $u$ is connected to $v$ via an edge}\\
&0,\quad\text{otherwise}.
\end{aligned}\right.
$$
It is clear to see that such kind of representation encapsules some basic information about underlying structure of a graph itself. For instance, the degree $k_{u}$ of vertex $u$ is equivalent to $k_{u}=\sum_{v=1}^{|\mathcal{V}|}a_{uv}$. To make further progress, the diagonal matrix, denoted by $\mathbf{D}$, may be viewed as follows: the $i$th diagonal entry is $k_{i}$, while all non-diagonal entries are zero, that is, $\mathbf{D}=\text{diag}[k_{1},k_{2},\dots,k_{|\mathcal{V}|}]$. Based on this, the corresponding Laplacian matrix, denoted by $\mathbf{L}$, is expressed as $\mathbf{L}=\mathbf{D}-\mathbf{A}$. Accordingly, the normalized version is $\mathbf{\mathcal{L}}=\mathbf{I}-\mathbf{D}^{-1}\mathbf{A}$ where $\mathbf{I}$ is an identity matrix with proper cardinality and $\mathbf{D}^{-1}$ indicates the inverse of matrix $\mathbf{D}$. For our propose, let $\lambda_{1},\lambda_{2},\cdots,\lambda_{|\mathcal{V}|}$ indicate the $|\mathcal{V}|$ eigenvalues of Laplacian matrix $\mathbf{L}$, which can be rearranged in an increasing order as follows, $0=\lambda_{1}<\lambda_{2}<\cdots<\lambda_{|\mathcal{V}|}$.

\subsection{Shortest path length}

Given a pair of vertices, say $u$ and $v$, in graph $\mathcal{G}(\mathcal{V},\mathcal{E})$, the distance $d_{uv}$, also usually thought of as shortest path length, is the edge number of any shortest path joining vertices $u$ and $v$. The diameter is the maximum among all the distances in graph. The sum over distances of all possible vertex pairs in graph $\mathcal{G}(\mathcal{V},\mathcal{E})$, denoted by $\mathcal{W}$ which is also called Wiener index, is by definition given by

\begin{equation}\label{eqa:MF-2-1-1}
\mathcal{W}=\frac{1}{2}\sum_{u\in \mathcal{V}}\sum_{v\in \mathcal{V}}d_{uv}.
\end{equation}
At meantime, we can define the mean shortest path length $\overline{\mathcal{W}}$ in the following form

\begin{equation}\label{eqa:MF-2-1-2}
\overline{\mathcal{W}}=\frac{\mathcal{W}}{|\mathcal{V}|(|\mathcal{V}|-1)/2}.
\end{equation}

\subsection{Random walks on graph}

When considering random walks on network, an important and fundamental parameter is first-passage time ($FPT$) for an arbitrarily given pair of vertices $u$ and $v$. By definition, the first-passage time from vertex $u$ to $v$, denoted by $\mathcal{F}_{u\rightarrow v}$, is the expected time taken by a walker starting out from vertex $u$ to first reach its destination vertex $v$. As before, for network $\mathcal{G}(\mathcal{V},\mathcal{E})$ as a whole, we can define two analogs relevant to $FPT$ for random walks, i.e.,

\begin{equation}\label{eqa:MF-2-2-1}
\mathcal{F}=\sum_{u}\sum_{v}\mathcal{F}_{u\rightarrow v},
\end{equation}
and,
\begin{equation}\label{eqa:MF-2-2-2}
\overline{\mathcal{F}}=\frac{\mathcal{F}}{|\mathcal{V}|(|\mathcal{V}|-1)}.
\end{equation}
It is worth mentioning that in general, the quantity $\mathcal{F}_{u\rightarrow v}$ is not necessarily equal to $\mathcal{F}_{v\rightarrow u}$ when a walker does random walks on network. Thus, the factor of $2$ is eliminated in the preceding two equations compared with Eqs.(\ref{eqa:MF-2-1-1}) and (\ref{eqa:MF-2-1-2}). The parameter $\overline{\mathcal{F}}$ is often called mean first-passage time ($MFPT$) for the sake of simplicity.

If network that we are discussing is a tree $\mathcal{T}$, there is an expression for $MFPT$ based on the eigenvalues of Laplacian matrix $\mathbf{L_{\mathcal{T}}}$, which is given by

\begin{equation}\label{eqa:MF-2-2-3}
\overline{\mathcal{F}}_{\mathcal{T}}=2\sum_{i=2}^{|\mathcal{V}_{\mathcal{T}}|}\frac{1}{\lambda_{i}}
\end{equation}
where symbol $|\mathcal{V}_{\mathcal{T}}|$ represents the total number of vertices in tree $\mathcal{T}$. As seen in the literature \cite{Vicsek-1983,Sheng-2019}, most of previous work concerning with quantity $MFPT$ on trees with intriguing properties, such as, fractal feature, is based on Eq.(\ref{eqa:MF-2-2-3}). Generally speaking, it is a standard technique for all trees to study mean first-passage time using Laplacian spectra.

\subsection{Electric network}

For a given graph $\mathcal{G}(\mathcal{V},\mathcal{E})$, we can construct a corresponding electrical network, simply referred to as $\mathcal{G}^{\dagger}(\mathcal{V}^{\dagger},\mathcal{E}^{\dagger})$, by replacing each edge in $\mathcal{E}$ with a unit resistor. For any two distinct vertices $u$ and $v$ in $\mathcal{V}^{\dagger}$, the effective resistance $\Omega_{uv}$ between them is defined as the
potential difference between $u$ and $v$ when a unit current from $u$ to $v$ is maintained. In the case $u$ identical to $v$, $\Omega_{uv}$ is equal to zero.

As known, effective resistance is in fact a measure of distance. So, the sum over effective resistances $\Omega_{uv}$ of all possible vertex pairs can be expressed as

\begin{equation}\label{eqa:MF-2-4-1}
\mathcal{R}_{\mathcal{G}^{\dagger}}=\sum_{u\in \mathcal{V}^{\dagger}}\sum_{v\in \mathcal{V}^{\dagger}}\Omega_{uv}.
\end{equation}
The quantity $\mathcal{R}_{\mathcal{G}^{\dagger}}$ is commonly called Kirchhoff index of electrical network $\mathcal{G}^{\dagger}(\mathcal{V}^{\dagger},\mathcal{E}^{\dagger})$. Analogously, Kirchhoff index $\mathcal{R}_{\mathcal{G}^{\dagger}}$ is also calculated in terms of the eigenvalues of Laplacian matrix $\mathbf{L_{\mathcal{G}}}$ of its underlying graph $\mathcal{G}(\mathcal{V},\mathcal{E})$, which is as below

\begin{equation}\label{eqa:MF-2-4-2}
\mathcal{R}_{\mathcal{G}^{\dagger}}=2|\mathcal{V}|\sum_{i=2}^{|\mathcal{V}|}\frac{1}{\lambda_{i}}.
\end{equation}

The mean effective resistance $\overline{\mathcal{R}_{\mathcal{G}^{\dagger}}}$ in $\mathcal{G}^{\dagger}(\mathcal{V}^{\dagger},\mathcal{E}^{\dagger})$, called network criticality, is read to

\begin{equation}\label{eqa:MF-2-4-2}
\overline{\mathcal{R}_{\mathcal{G}^{\dagger}}}=\frac{\mathcal{R}_{\mathcal{G}^{\dagger}}}{|\mathcal{V}|(|\mathcal{V}|-1)}=\frac{2}{|\mathcal{V}|-1}\sum_{i=2}^{|\mathcal{V}|}\frac{1}{\lambda_{i}}.
\end{equation}
$\overline{\mathcal{R}_{\mathcal{G}^{\dagger}}}$ quantifies the network robustness of $\mathcal{G}(\mathcal{V},\mathcal{E})$ as a communication network: smaller value for $\overline{\mathcal{R}_{\mathcal{G}^{\dagger}}}$ implies that network $\mathcal{G}(\mathcal{V},\mathcal{E})$ is more robust. For brevity, we make use of the unique symbol $\mathcal{G}(\mathcal{V},\mathcal{E})$ to represent a graph and its corresponding electrical network in the remainder of this paper.

\section{Graphic operation and Framework}

Here, we propose some graphic operations used later. It is worth noting that the below operations can still be generalized in many manners. Some simple generalized versions are shown in the rest of this section. Nonetheless, we aim at clarifying the thought behind these operations from the perspective of methodology. That is to say, it suffices to only consider graphic operations defined below. For ease of exposition, let us first define four vectors, say $\overrightarrow{m}_{\mu}$, $\overrightarrow{n}_{\nu}$, $\overrightarrow{p}_{\mu}$ and $\overrightarrow{q}_{\nu}$, in the following form

$$\overrightarrow{m}_{\mu}=(m_{1},m_{2},\cdots,m_{\mu}),\qquad \overrightarrow{n}_{\nu}=(n_{1},n_{2},\cdots,n_{\nu}),\qquad \overrightarrow{p}_{\mu}=(p_{1},p_{2},\cdots,p_{\mu}),\qquad \overrightarrow{q}_{\nu}=(q_{1},q_{2},\cdots,q_{\nu}).$$
Additionally, we require that the last two vectors, $\overrightarrow{p}$ and $\overrightarrow{q}$, meet the following criteria

$$\forall i \quad p_{i}\geq0,\quad \text{and}\quad ||\overrightarrow{p}_{\mu}||_{1}=1;\qquad \forall j \quad q_{j}\geq0,\quad \text{and}\quad ||\overrightarrow{q}_{\nu}||_{1}=1$$
where $||\bullet||_{1}$ is the $\ell_{1}$-norm. It is clear to see that vectors $\overrightarrow{p}_{\mu}$ and $\overrightarrow{q}_{\nu}$ are in the $\mu$-dimensional probability simplex and the $\nu$-dimensional probability simplex, respectively. Next, each entry $m_{i}$ in vector $\overrightarrow{m}_{\mu}$ is an arbitrary natural number at least larger than $0$. At the same time, we demand $m_{i}\neq m_{j}$ for an arbitrary pair of distinct indices $i$ and $j$. The similar requirement holds on vector $\overrightarrow{n}_{\nu}$.

\begin{figure}
\centering
  \includegraphics[height=3cm]{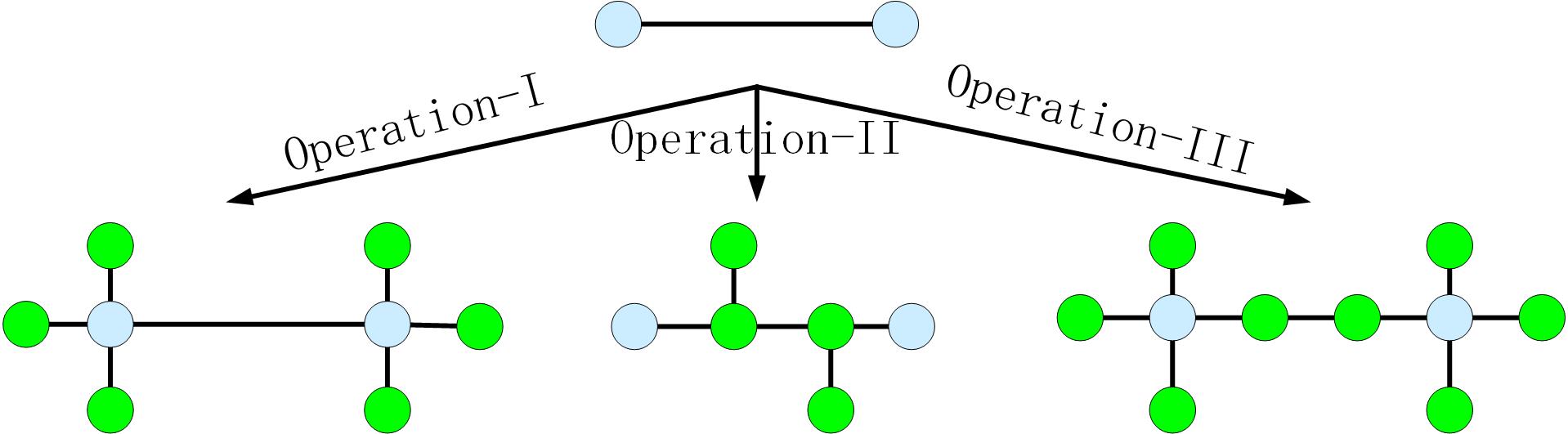}\\
  \vskip0.5cm
{\small Fig.1. (Color online) The diagrams of three operations. We only show deterministic versions corresponding to three operations for convenience. Particularly, the seed is an edge $uv$. When executing Operation-I, we assume that $\mu=3$ and $m_{i}=1$ for all $i\in[1,3]$ as shown in the leftmost panel. Similarly, assume that $\mu=\nu=1$, $m_{1}=2$ and $n_{1}=1$ in Operation-II, the intermediate panel shows an example. In the third operation as shown in the rightmost panel, we make use of $\mu=4$ and $m_{i}=1$ for all $i\in[1,4]$. Note that red edges represent those newly created edges during implementing operation.}
\end{figure}

\textbf{Definition 1} For a given graph $\mathcal{G}(\mathcal{V},\mathcal{E})$, \emph{Operation-I} is to attach $\mu$ paths to each vertex $v$, each having length $m_{i}-1$ with probability $p_{i}$. This operation is called \emph{Vertex-based uniform growth mechanism} ($VUGM$) mainly because such an operation is applied on each vertex without considering the structural properties of vertex itself, such as, vertex degree.

\emph{Example 1} When the seed is an edge and each path degenerates into an isolated vertex, the resulting graph is in fact the deterministic uniform growth tree, denoted by $\mathcal{Y}_{I}(t)$, after iteratively manipulating Operation-I $t$ times. This kind of trees have been well studied due to their own prevalent applications in real-life world \cite{Masuda-2017}, for instance, modeling epidemic spread.

\textbf{Definition 2} Given a graph $\mathcal{G}(\mathcal{V},\mathcal{E})$, we can insert all centres of $m_{i}$ star-like graphs\footnote[1]{A star-like graph $\mathcal{G}(\mathcal{V},\mathcal{E})$ is one graph that shares a similar topological structure to star graph \cite{Bondy-2008}. That is to say, there exists a vertex acting as the center in star-like graph $\mathcal{G}(\mathcal{V},\mathcal{E})$, and the central vertex is attached to some paths. Vividly, each path is also viewed as a ``tentacle" for brevity. An illustrative example is plotted in Appendix.} into each edge $uv$ with probability $p_{i}$ where each star-like graph has $\nu$ ``tentacles". Here, each ``tentacle" owns $n_{j}$ vertices with probability $q_{j}$ independently. Such a manipulation is defined as \emph{Operation-II}, which is also viewed as \emph{Edge-based uniform growth mechanism} ($EUGM$).

\emph{Example 2} Similarly, if the seed is an edge, and $m_{i}$ and $n_{j}$ are equal to $1$ for all $i$ and $j$, we can obtain the well-known $T$-graph, denoted by $\mathcal{Y}_{II}(t)$, through iteratively executing Operation-II $t$ times. Some structural properties on $T$-graph have been discussed in \cite{Kahng-1989,Agliari-2008} particularly because it may serve as a simple model illustrating the inhomogeneity and scale-invariance of many disordered materials in physics.

\textbf{Definition 3} Considering a graph $\mathcal{G}(\mathcal{V},\mathcal{E})$ with largest vertex degree $k_{max}$ which is no greater than parameter $\mu$, we first insert $m_{i}$ new vertices into each edge incident with vertex $v$ according to probability $p_{i}$. Next, two cases need to be considered: (1) If the degree $k_{v}$ of vertex $v$ is equal to $\mu$, then there is nothing more to do; (2) On the other hand, i.e., $k_{v}<\mu$, we need to connect $\mu-k_{v}$ paths to vertex $v$ in which each newly added path has $m_{i}$ vertices in terms of probability $p_{i}$. The procedure above is referred to as \emph{Operation-III}, which is called \emph{mixture uniform growth mechanism} ($MUGM$) as well.

\emph{Example 3} We select a star $S_{\mu}$ on $\mu+1$ vertices as a seed, assume that each entry $m_{i}$ in vector $\overrightarrow{m}_{\mu}$ is equivalent to $1$, and obtain the famous Vicsek fractal $V_{1}^{\mu}(t)$ after running Operation-III $t$ times. The key parameters pertaining to Vicsek fractal $V_{1}^{\mu}(t)$, such as, mean first-passage time, have been widely studied in the published papers including Ref.\cite{Zhang-2010}. One of most important reasons for this is that Vicsek fractal $V_{1}^{\mu}(t)$ can be utilized to describe the underlying structure of some polymers in chemistry.

Figure 1 shows some examples in order to better understand details about operations introduced above. As mentioned above, the goal of this paper is to study stochastic uniform growth trees. Hence, the initial graph (aka seed) is always an arbitrary tree $\mathcal{T}$. In Examples 1-3, the selected seeds are in fact trees with specific structural property. In the following, we are going to establish a principled framework $\Upsilon$ based on three different types of operations stated in Defs.1-3.

\textbf{Framework $\Upsilon$}

At $t=0$, an arbitrary tree $\mathcal{T}$ is chosen as the seed and four vectors are defined as above.

At $t=1$, one of three operations built is applied on tree $\mathcal{T}$. The resulting tree is denoted by $\mathcal{T}(1)$. More specifically, applying Operation-I yields tree $\mathcal{T}_{I}(1)$ and similarly for other two operations.

At $t\geq2$, the new generation $\mathcal{T}(t)$ can be obtained from the preceding tree $\mathcal{T}(t-1)$ by implementing the same operation as in the previous time step.

We illustrate the principled framework $\Upsilon$ in Fig.2 for the goal of expounding concrete procedures. Specifically speaking, this framework outputs three distinct families of stochastic uniform growth trees, that is,  $\mathcal{T}_{I}(t)$,  $\mathcal{T}_{II}(t)$ as well as $\mathcal{T}_{III}(t)$. In principle, we are able to replace each added path by an arbitrary tree when applying Operation-I. This leads to more general stochastic uniform growth trees. Analogously, some generalization technologies can also be adopted in other two operations for creating generalized versions. Due to space limitation, we omit it here. It is worth noticing that the goal of this work is to provide a guideline for creating stochastic uniform growth trees. In what follows, let us pay considerable attention on studying structural parameters on the resulting tree networks.

\begin{figure}
\centering
  \includegraphics[height=4cm]{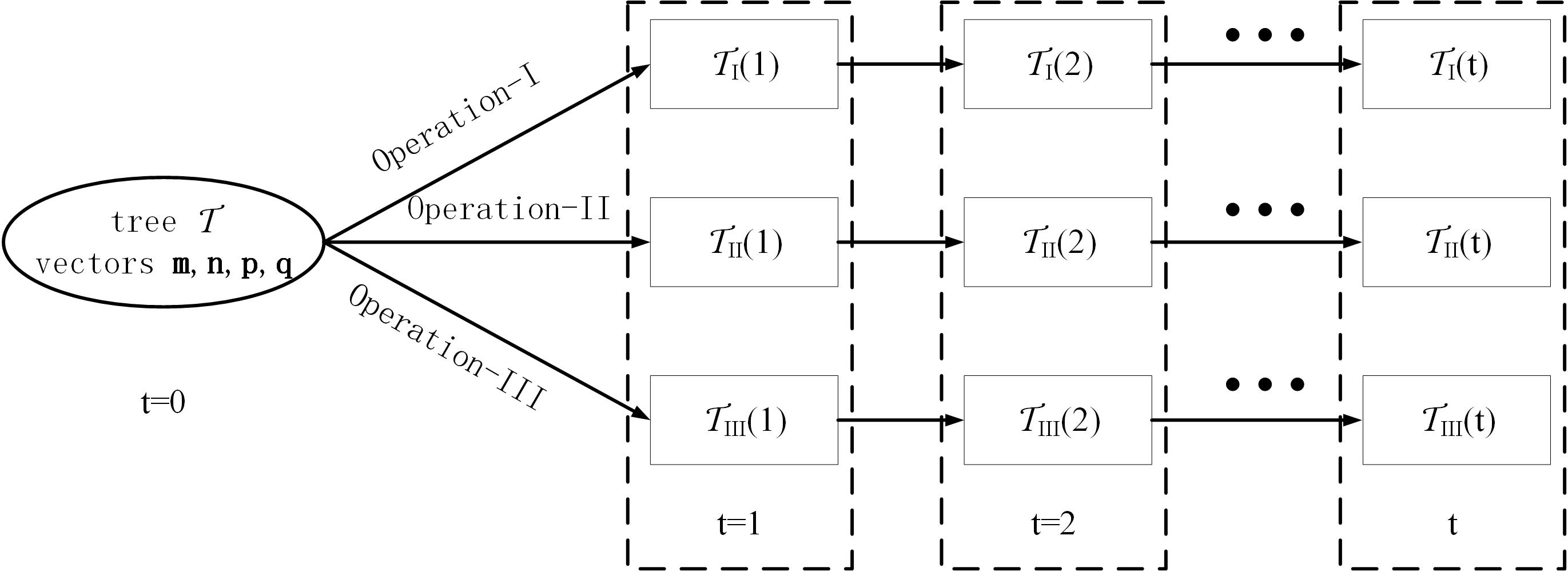}\\
  \vskip0.5cm
{\small Fig.2. (Color online) The diagram of the principled framework $\Upsilon$ for creating stochastic uniform growth trees.}
\end{figure}

\section{Main Results}

The goal of this section is to display main results. For instance, we determine the analytic solutions to mean first-passage time for random walks on three different kinds of stochastic uniform growth trees output by framework $\Upsilon$. As an immediate consequence, those previously published formulas for deterministic versions can be easily obtained by substituting available parameters into our results.

First of all, let us take a lemma.

\textbf{Lemma 1 \cite{Ma-2020-1}} For a designated vertex $v$ in graph $\mathcal{G}(\mathcal{V},\mathcal{E})$, the sum $\mathcal{F}_{v}$ of first-passage times for all vertices $u$ in its neighboring set $\mathcal{N}_{v}$ follows

\begin{equation}\label{Section-3-1-0}
\mathcal{F}_{v}=\sum_{u\in\mathcal{N}_{v}}\mathcal{F}_{u\rightarrow v}=2|\mathcal{E}|-k_{v}
\end{equation}
where $k_{v}$ represents the degree of vertex $v$ (see Eq.(14) in \cite{Ma-2020-1} for more details).

From Eq.(\ref{Section-3-1-0}), we can obtain a more helpful result for tree $\mathcal{T}$ as shown in the next proposition, which enables us to establish a proof of theorem 4.

\textbf{Proposition 2} Given an edge $uv$ in tree $\mathcal{T}$, the first-passage times $\mathcal{F}_{u\rightarrow v}$ and $\mathcal{F}_{v\rightarrow u}$ are given by

\begin{equation}\label{Section-3-2-0}
\mathcal{F}_{u\rightarrow v}=2|\mathcal{E}_{u}|+1,\qquad \mathcal{F}_{v\rightarrow u}=2|\mathcal{E}_{v}|+1
\end{equation}
where $|\mathcal{E}_{v}|$ is the total number of edges in the $v$-root tree and $|\mathcal{E}_{u}|$ has the same meaning for the $u$-root tree. Deleting edge $uv$ in tree $\mathcal{T}$ yields two subtrees. The subtree containing vertex $v$ is called $v$-root tree, and the other is viewed as  $u$-root tree.

\textbf{\emph{Proof}} The correctness of Eq.(\ref{Section-3-2-0}) is validated by deduction on edge number. First of all, proposition 2 is trivial when tree $\mathcal{T}$ is an edge. With loss of generality, we assume that the edge number in tree $\mathcal{T}$ is no less than $1$. It is clear to see that proposition 2 holds when either of two vertices $u$ and $v$ is a leaf vertex. For the sake of simplicity, we only prove the first equation. The other can be checked in a similar manner. Assume now that the first equation is true for root tree having edge number less than $|\mathcal{E}_{u}|$. By definition, the quantity $\mathcal{F}_{u\rightarrow v}$ is expressed as

\begin{equation}\label{Section-3-2-1}
\mathcal{F}_{u\rightarrow v}=\frac{1}{k_{u}}+\frac{1}{k_{u}}\sum_{u_{i}(\neq v)\in \mathcal{N}_{u}}(\mathcal{F}_{u_{i}\rightarrow v}+1).
\end{equation}
In view of nature of tree itself, Eq.(\ref{Section-3-2-1}) can be rewritten as

\begin{equation}\label{Section-3-2-2}
\mathcal{F}_{u\rightarrow v}=\frac{1}{k_{u}}+\frac{1}{k_{u}}\sum_{u_{i}(\neq v)\in \mathcal{N}_{u}}(\mathcal{F}_{u_{i}\rightarrow u}+\mathcal{F}_{u\rightarrow v}+1).
\end{equation}
Obviously, the edge number $|\mathcal{E}_{u_{i}}|$ of $u_{i}$-root tree is strictly less than quantity $|\mathcal{E}_{u}|$. By hypothesis, we can have
\begin{equation}\label{Section-3-2-3}
\mathcal{F}_{u\rightarrow v}=\frac{1}{k_{u}}+\frac{1}{k_{u}}\sum_{u_{i}(\neq v)\in \mathcal{N}_{u}}\left[(2|\mathcal{E}_{u_{i}}|+1)+\mathcal{F}_{u\rightarrow v}+1\right].
\end{equation}
Multiplying $k_{u}$ on the both-hand sides of Eq.(\ref{Section-3-2-3}) yields

\begin{equation}\label{Section-3-2-4}
k_{u}\mathcal{F}_{u\rightarrow v}=1+\sum_{u_{i}(\neq v)\in \mathcal{N}_{u}}\left[(2|\mathcal{E}_{u_{i}}|+1)+\mathcal{F}_{u\rightarrow v}+1\right].
\end{equation}
Next, we have

\begin{equation}\label{Section-3-2-5}
\begin{aligned}\mathcal{F}_{u\rightarrow v}&=k_{u}+\sum_{u_{i}(\neq v)\in \mathcal{N}_{u}}(2|\mathcal{E}_{u_{i}}|+1)\\
&=2\sum_{u_{i}(\neq v)\in \mathcal{N}_{u}}(|\mathcal{E}_{u_{i}}|+1)+1
\end{aligned}.
\end{equation}
According to nature of tree, we complete the proof of Eq.(\ref{Section-3-2-0}) upon Eq.(\ref{Section-3-2-5}).

\textbf{Corollary 3} For an arbitrarily given pair of vertices $u$ and $v$ in tree $\mathcal{T}$, the commute time  $\mathcal{F}_{v\leftrightarrow u}$ is given by

\begin{equation}\label{Section-3-3-0}
\mathcal{F}_{v\leftrightarrow u}=2|\mathcal{E}_{\mathcal{T}}|d_{uv}
\end{equation}
in which $d_{uv}$ is the distance between them. This is an immediate result of proposition 2 and we thus omit its proof here. Note that a more general version relevant to commute time $\mathcal{F}_{v\leftrightarrow u}$ on graph has been derived using spectral technique \cite{Guex-2015}. To make further progress, based on Eq.(\ref{Section-3-3-0}), we can establish a connection of Wiener index $\mathcal{W}_{\mathcal{T}}$ to mean first-passage time $\overline{\mathcal{F}}_{\mathcal{T}}$ as shown in the following theorem.

\textbf{Theorem 4} Consider random walks on a tree $\mathcal{T}$, there is a formula between two structural parameters $\mathcal{W}_{\mathcal{T}}$ and $\overline{\mathcal{F}}_{\mathcal{T}}$, as follows

\begin{equation}\label{Section-3-4-0}
2\mathcal{W}_{\mathcal{T}}=|\mathcal{V}_{\mathcal{T}}|\overline{\mathcal{F}}_{\mathcal{T}}.
\end{equation}
According to the simplicity of proof, we omit it here. Note also that a more general consequence corresponding to Eq.(\ref{Section-3-4-0}) is found in terms of spectral technique \cite{Chandra-1989} (see theorem 2.1 in \cite{Chandra-1989} for more details). Obviously, this implies that we have the ability to derive the analytic solutions to mean first-passage time for random walks on all the stochastic uniform growth trees if the derivation of Wiener index on corresponding trees is easy to deal with. We make a statement in advance that for convenience and brevity, we will select a tree $\mathcal{T}$ on $h$ vertices as a seed to create the candidates through framework $\Upsilon$ in the rest of this paper, and denote by $\mathcal{W}_{\mathcal{T}}$ the corresponding Wiener index of tree $\mathcal{T}$.

In what follows, we study three distinct kinds of stochastic uniform growth trees, i.e., $\mathcal{T}_{I}(t)$, $\mathcal{T}_{II}(t)$ and $\mathcal{T}_{III}(t)$, and derive the analytic solutions to some structural parameters including Wiener index and mean first-passage time. It should be mentioned that in view of randomness of trees $\mathcal{T}_{I}(t)$, $\mathcal{T}_{II}(t)$ and $\mathcal{T}_{III}(t)$, the results obtained are expected expressions. On the other hand, we attempt to make use of a brief yet unambiguous presentation in the following discussions. For instance, we use ``solution of Wiener index" instead of ``expected solution of Wiener index".

\subsection{Tree $\mathcal{T}_{I}(t)$}

\textbf{Theorem 5} The solution of Wiener index $\mathcal{W}_{\mathcal{T}_{I}(1)}$ of tree $\mathcal{T}_{I}(1)$ is given by

\begin{equation}\label{Section-4-1-0}
\begin{aligned}\mathcal{W}_{\mathcal{T}_{I}(1)}&=\mathcal{W}_{\mathcal{T}}\left(\mu\sum_{i=1}^{\mu}p_{i}m_{i}+1\right)^{2}+\mu\left(\mu\sum_{i=1}^{\mu}p_{i}m_{i}+1\right)\sum_{i=1}^{\mu}p_{i}\left(
                                                             \begin{array}{c}
                                                               m_{i}+1 \\
                                                               2 \\
                                                             \end{array}
                                                           \right)h^{2}\\
&\quad+\left\{\mu\sum_{i}^{\mu}p_{i}\sum_{l=2}^{m_{i}}\left(
                                                             \begin{array}{c}
                                                               l \\
                                                               2 \\
                                                             \end{array}
                                                           \right)+\left[\sum_{i=1}^{\mu}2p_{i}m_{i}\left(
                                                             \begin{array}{c}
                                                               \mu \\
                                                               2 \\
                                                             \end{array}
                                                           \right)-\mu^{2}\sum_{i=1}^{\mu}p_{i}m_{i}\right]\sum_{i=1}^{\mu}p_{i}\left(
                                                             \begin{array}{c}
                                                               m_{i}+1 \\
                                                               2 \\
                                                             \end{array}
                                                           \right)\right\}h
\end{aligned}.
\end{equation}
For convenience, the formula above is commonly referred to as the $\mathcal{W}$-polynomial for Wiener index of tree $\mathcal{T}_{I}(1)$ whose variables are parameters $\mathcal{W}_{\mathcal{T}}$ and $h$ of seed $\mathcal{T}$. Roughly speaking, such a polynomial is able to be expressed as

\begin{equation}\label{Section-4-1-0-0}
f_{I}(\mathcal{W}_{\mathcal{T}},h)\triangleq a_{I}\mathcal{W}_{\mathcal{T}}+b_{I}h^{2}+c_{I}h.
\end{equation}

\textbf{\emph{Proof}} Before beginning with our discussions, some necessary notations are introduced as follows. We denote by $\Lambda_{u}^{I}$ the set of vertices which are added into tree $\mathcal{T}_{I}(1)$ by applying $VUGM$ to vertex $u$ in seed $\mathcal{T}$ \footnote[2]{As an example, $\Lambda_{u}^{I}$ is composed of threes green vertices connected to vertex $u$ in the left-most panel of Fig.1.}. This suggests that all the newly created vertices belong to set $\bigcup_{u\in\mathcal{T}}\Lambda_{u}^{I}$, and then the set $V_{\mathcal{T}_{I}(1)}/\bigcup_{u\in\mathcal{T}}\Lambda_{u}^{I}$ contains all the vertices previously belonging to tree $\mathcal{T}$. For simplicity of presentation, we define $\Omega_{1}^{I}=V_{\mathcal{T}_{I}(1)}/\bigcup_{u\in\mathcal{T}}\Lambda_{u}^{I}$. Based on this, the concrete demonstrations can be shown in stages.

\emph{Case 1} For an arbitrary pair of vertices, say $u$ and $v$, in vertex set $\Omega_{1}^{I}$, the distance between them $d'_{uv}$ keeps unchanged after applying $VUGM$ to seed $\mathcal{T}$. Thus, it is straightforward to see

\begin{equation}\label{Section-4-1-0-1}
\mathcal{W}_{\mathcal{T}_{I}(1)}(1)\triangleq\frac{1}{2}\sum_{u\in\Omega_{1}^{I}}\sum_{v\in\Omega_{1}^{I}}d'_{uv}=\mathcal{W}_{\mathcal{T}}.
\end{equation}

\emph{Case 2} Similarly, we can without difficulty obtain the following formula

\begin{equation}\label{Section-4-1-0-2}
\begin{aligned}\mathcal{W}_{\mathcal{T}_{I}(1)}(2)&\triangleq\sum_{u\in\Omega_{1}^{I}}\sum_{u_{\alpha}\in \Lambda^{I}_{u}}d'_{uu_{\alpha}}+\frac{1}{2}\sum_{u\in\Omega_{1}^{I}}\sum_{u_{\alpha}\in \Lambda^{I}_{u}}\sum_{u_{\beta}\in \Lambda^{I}_{u}}d'_{u_{\alpha}u_{\beta}}\\
&=h\left[\mu\sum_{i=1}^{\mu}p_{i}\left(
                                                             \begin{array}{c}
                                                               m_{i}+1 \\
                                                               2 \\
                                                             \end{array}
                                                           \right)+\sum_{i=1}^{\mu}2p_{i}m_{i}\left(
                                                             \begin{array}{c}
                                                               \mu \\
                                                               2 \\
                                                             \end{array}
                                                           \right)\sum_{i=1}^{\mu}p_{i}\left(
                                                             \begin{array}{c}
                                                               m_{i}+1 \\
                                                               2 \\
                                                             \end{array}
                                                           \right)+\mu\sum_{i}^{\mu}p_{i}\sum_{l=2}^{m_{i}}\left(
                                                             \begin{array}{c}
                                                               l \\
                                                               2 \\
                                                             \end{array}
                                                           \right)\right]
                                                           \end{aligned},
\end{equation}
in which the summation over distances between arbitrary pair of vertices $u_{\alpha}$ and $u_{\beta}$ in set $\Lambda^{I}_{u}$ is viewed as
$$\frac{1}{2}\sum_{u_{\alpha}\in \Lambda^{I}_{u}}\sum_{u_{\beta}\in \Lambda^{I}_{u}}d'_{u_{\alpha}u_{\beta}}.$$
At the same time, it is worth noticing that for $a<b$, we assume $\left(
                                                             \begin{array}{c}
                                                               a \\
                                                               b \\
                                                             \end{array}
                                                           \right)=0$ and $\sum_{x=b}^{a}x=0.$

\emph{Case 3} Now, let us focus on the derivation of sum over distances of all possible vertex pairs $u$ in $\Omega_{1}^{I}$ and $v_{\alpha}$ in $\Lambda^{I}_{v}$ where $u$ is distinct from $v$. Without loss of generality, suppose that the path joining vertex $u$ to $v_{\alpha}$ in tree $\mathcal{T}_{I}(1)$ is defined as $\mathcal{P}_{uv_{\alpha}}$. Consequently, it is clear to find path $\mathcal{P}_{uv_{\alpha}}$ to contain two sub-paths, say $\mathcal{P}_{uv}$ and $\mathcal{P}_{vv_{\alpha}}$. Thus, there is an identity $d'_{uv_{\alpha}}=d'_{uv}+d'_{vv_{\alpha}}.$ Accordingly, we can obtain

\begin{equation}\label{Section-4-1-0-3}
\begin{aligned}\mathcal{W}_{\mathcal{T}_{I}(1)}(3)&\triangleq\sum_{u\in\Omega_{1}^{I}}\sum_{v(\neq u)\in\Omega_{1}^{I}}\sum_{v_{\alpha}\in \Lambda^{I}_{v}}d'_{uv_{\alpha}}\\
&=2\mu\left[\left(
                                                             \begin{array}{c}
                                                               h \\
                                                               2 \\
                                                             \end{array}
                                                           \right)\sum_{i=1}^{\mu}p_{i}\left(
                                                             \begin{array}{c}
                                                               m_{i}+1 \\
                                                               2 \\
                                                             \end{array}
                                                           \right)+\sum_{i=1}^{\mu}p_{i}m_{i}\mathcal{W}_{\mathcal{T}}\right]
                                                           \end{aligned}.
\end{equation}

\emph{Case 4} The remainder of this task is to calculate the distance between vertices $u_{\alpha}$ and $v_{\beta}$, denoted by $d'_{u_{\alpha}v_{\beta}}$. Note that this pair of vertices are from two different sets $\Lambda^{I}_{u}$ and $\Lambda^{I}_{v}$, respectively. By analogy with calculation of quantity $\mathcal{W}_{\mathcal{T}_{I}(1)}(3)$ in the previous case, the sum over distances of this type is calculated to yield

\begin{equation}\label{Section-4-1-0-4}
\begin{aligned}\mathcal{W}_{\mathcal{T}_{I}(1)}(4)&\triangleq\sum_{u\in\Omega_{1}^{I}}\sum_{v(\neq u)\in\Omega_{1}^{I}}\sum_{u_{\alpha}\in\Lambda^{I}_{u}}\sum_{v_{\beta}\in \Lambda^{I}_{v}}d'_{u_{\alpha}v_{\beta}}\\
&=\left(\mu\sum_{i=1}^{\mu}p_{i}m_{i}\right)^{2}\mathcal{W}_{\mathcal{T}}+2\mu\sum_{i=1}^{\mu}p_{i}m_{i}\left(
                                                             \begin{array}{c}
                                                               h \\
                                                               2 \\
                                                             \end{array}
                                                           \right)\left[\mu\sum_{i=1}^{\mu}p_{i}\left(
                                                             \begin{array}{c}
                                                               m_{i}+1 \\
                                                               2 \\
                                                             \end{array}
                                                           \right)\right]
\end{aligned}.
\end{equation}

Armed with all the cases, we prove Eq.(\ref{Section-4-1-0}) based on summation $\mathcal{W}_{\mathcal{T}_{I}(1)}=\sum_{i=1}^{4}\mathcal{W}_{\mathcal{T}_{I}(1)}(i)$ after performing some fundamental arithmetics. $\hfill\Box$

In fact, theorem 5 provides us with an approach to determining the analytic solution to Wiener index $\mathcal{W}_{\mathcal{T}_{I}(t)}$ on stochastic uniform growth tree $\mathcal{T}_{I}(t)$. Now, the only requirement is to first know vertex number $\mathcal{V}_{\mathcal{T}_{I}(t)}$ of tree $\mathcal{T}_{I}(t)$. This is easily derived via the following recursive relation

$$|\mathcal{V}_{\mathcal{T}_{I}(t)}|=\left(1+\mu\sum_{i=1}^{\mu}p_{i}m_{i}\right)|\mathcal{V}_{\mathcal{T}_{I}(t-1)}|.$$
Using the initial condition $|\mathcal{V}_{\mathcal{T}_{I}(0)}|=|\mathcal{V}_{\mathcal{T}}|=h$, we obtain

\begin{equation}\label{Section-4-1-1}
|\mathcal{V}_{\mathcal{T}_{I}(t)}|=h\left(1+\mu\sum_{i=1}^{\mu}p_{i}m_{i}\right)^{t}.
\end{equation}
With the results mentioned above, we reach the next proposition.

\textbf{Proposition 6} The solution of Wiener index $\mathcal{W}_{\mathcal{T}_{I}(t)}$ of tree $\mathcal{T}_{I}(t)$ is given by

\begin{equation}\label{Section-4-1-2}
\begin{aligned}\mathcal{W}_{\mathcal{T}_{I}(t)}&=\mathcal{W}_{\mathcal{T}}\left(\mu\sum_{i=1}^{\mu}p_{i}m_{i}+1\right)^{2t}+\mu t h^{2}\sum_{i=1}^{\mu}p_{i}\left(
                                                             \begin{array}{c}
                                                               m_{i}+1 \\
                                                               2 \\
                                                             \end{array}
                                                           \right)\left(1+\mu\sum_{i=1}^{\mu}p_{i}m_{i}\right)^{2t-1}\\
                                                          &+h\left\{\sum_{i}^{\mu}p_{i}\sum_{l=2}^{m_{i}}\left(
                                                             \begin{array}{c}
                                                               l \\
                                                               2 \\
                                                             \end{array}
                                                           \right)+\left[\sum_{i=1}^{\mu}2p_{i}m_{i}\left(
                                                             \begin{array}{c}
                                                               \mu \\
                                                               2 \\
                                                             \end{array}
                                                           \right)-\mu^{2}\sum_{i=1}^{\mu}p_{i}m_{i}\right]\sum_{i=1}^{\mu}p_{i}\left(
                                                             \begin{array}{c}
                                                               m_{i}+1 \\
                                                               2 \\
                                                             \end{array}
                                                           \right)\right\}\Phi^{I}_{1}
\end{aligned},
\end{equation}
in which the solution of symbol $\Phi^{I}_{1}$ is

$$\Phi^{I}_{1}=\frac{\left(1+\mu\sum_{i=1}^{\mu}p_{i}m_{i}\right)^{2t-1}-\left(1+\mu\sum_{i=1}^{\mu}p_{i}m_{i}\right)^{t-1}}{\mu\sum_{i=1}^{\mu}p_{i}m_{i}}.$$

Based on Eqs.(\ref{Section-4-1-0}) and (\ref{Section-4-1-1}), the aforementioned proposition can be proved in an iterative manner, and thus we omit it here. Using the relation shown in Eq.(\ref{Section-3-4-0}), the closed-form solution to quantity $\overline{\mathcal{F}}_{\mathcal{T}_{I}(t)}$ is obtained with respect to Eq.(\ref{Section-4-1-2}) immediately, which is shown in the next theorem.

\textbf{Theorem 7} The analytic solution to mean first-passage time $\overline{\mathcal{F}}_{\mathcal{T}_{I}(t)}$ for random walks on tree $\mathcal{T}_{I}(t)$ is given by

\begin{equation}\label{Section-4-1-3}
\begin{aligned}\overline{\mathcal{F}}_{\mathcal{T}_{I}(t)}&=\frac{2\mathcal{W}_{\mathcal{T}}}{h}\left(\mu\sum_{i=1}^{\mu}p_{i}m_{i}+1\right)^{t}+ 2\mu th\sum_{i=1}^{\mu}p_{i}\left(
                                                             \begin{array}{c}
                                                               m_{i}+1 \\
                                                               2 \\
                                                             \end{array}
                                                           \right)\left(1+\mu\sum_{i=1}^{\mu}p_{i}m_{i}\right)^{t-1}\\
&+2h\left\{\sum_{i}^{\mu}p_{i}\sum_{l=2}^{m_{i}}\left(
                                                             \begin{array}{c}
                                                               l \\
                                                               2 \\
                                                             \end{array}
                                                           \right)+\left[\sum_{i=1}^{\mu}2p_{i}m_{i}\left(
                                                             \begin{array}{c}
                                                               \mu \\
                                                               2 \\
                                                             \end{array}
                                                           \right)-\mu^{2}\sum_{i=1}^{\mu}p_{i}m_{i}\right]\sum_{i=1}^{\mu}p_{i}\left(
                                                             \begin{array}{c}
                                                               m_{i}+1 \\
                                                               2 \\
                                                             \end{array}
                                                           \right)\right\}\Phi^{I}_{2}
\end{aligned},
\end{equation}
where we have made use of

$$\Phi^{I}_{2}=\frac{\left(1+\mu\sum_{i=1}^{\mu}p_{i}m_{i}\right)^{t-1}-\left(1+\mu\sum_{i=1}^{\mu}p_{i}m_{i}\right)^{-1}}{\mu\sum_{i=1}^{\mu}p_{i}m_{i}}.$$

This is now an obvious consequence and hence we omit its proof. In the limit of large graph size, the scaling relation between two structural parameters, $V_{\mathcal{T}_{I}(t)}$ and $\overline{\mathcal{F}}_{\mathcal{T}_{I}(t)}$, obeys

\begin{equation}\label{Section-4-1-4}
\overline{\mathcal{F}}_{\mathcal{T}_{I}(t)}=O\left(\Theta_{I} |\mathcal{V}_{\mathcal{T}_{I}(t)}|\right),\qquad\qquad \Theta_{I}=\frac{\mu t\sum_{i=1}^{\mu}p_{i}m_{i}(m_{i}-1)}{1+\mu\sum_{i=1}^{\mu}p_{i}m_{i}}.
\end{equation}
This implies that the mean first-passage time for random walks on stochastic uniform growth tree $\mathcal{T}_{I}(t)$ is approximately linearly correlated with the total number of vertices as $t\rightarrow \infty$.

As mentioned above, a special member in tree family $\mathcal{T}_{I}(t)$, namely, deterministic version $\mathcal{Y}_{I}(t)$, has been studied analytically. Now, by using our formula in Eq.(\ref{Section-4-1-3}), the corresponding theoretical expression of mean first-passage time is timely derived as shown in corollary 8.

\textbf{Corollary 8} The exact solution to mean first-passage time $\overline{\mathcal{F}}_{\mathcal{Y}_{I}(t)}$ on deterministic uniform growth tree $\mathcal{Y}_{I}(t)$ is

\begin{equation}\label{Section-4-1-5}
\overline{\mathcal{F}}_{\mathcal{Y}_{I}(t)}=(4\mu t+\mu-1)(1+\mu)^{t-1}+\frac{2}{1+\mu}.
\end{equation}
This is completely the same as the previous result in published paper \cite{Wang-2012} (see Eq.(31) in \cite{Wang-2012} for more details), suggesting that theorem 7 is sound.

\subsection{Tree $\mathcal{T}_{II}(t)$}

\textbf{Theorem 9} The solution of Wiener index $\mathcal{W}_{\mathcal{T}_{II}(1)}$ of tree $\mathcal{T}_{II}(1)$ is given by

\begin{equation}\label{Section-4-2-0}
\begin{aligned}\mathcal{W}_{\mathcal{T}_{II}(1)}&=\mathcal{W}_{\mathcal{T}}\Phi_{4}^{II}\left(\Phi_{1}^{II}+1\right)^{2}+\left[\Phi_{1}^{II}\Phi_{3}^{II}-\Phi_{4}^{II}\left(\Phi_{1}^{II}\right)^{2}-\Phi_{4}^{II}\Phi_{1}^{II}+\Phi_{3}^{II}\right]h^{2}\\
&\quad+\left[2\Phi_{4}^{II}\left(\Phi_{1}^{II}\right)^{2}-3\Phi_{1}^{II}\Phi_{3}^{II}+\Phi_{4}^{II}\Phi_{1}^{II}+\Phi_{2}^{II}-\Phi_{3}^{II}\right]h+\left[2\Phi_{1}^{II}\Phi_{3}^{II}-\Phi_{4}^{II}\left(\Phi_{1}^{II}\right)^{2}-\Phi_{2}^{II}\right]
\end{aligned}.
\end{equation}
Note that the concrete meanings of symbols $\Phi_{i}^{II}$ ($i\in[1,4]$) are deferred in the next proof for the sake of argument. As mentioned previously, the $\mathcal{W}$-polynomial for tree $\mathcal{W}_{\mathcal{T}_{II}(1)}$ is written as

\begin{equation}\label{Section-4-2-0-0}
f_{II}(\mathcal{W}_{\mathcal{T}},h)\triangleq a_{II}\mathcal{W}_{\mathcal{T}}+b_{II}h^{2}+c_{II}h+d_{II}.
\end{equation}

\textbf{\emph{Proof}} Let us first recall $EUGM$, and then find that at present, operation is applied to each edge in seed $\mathcal{T}$. For our purpose, we use $\Omega_{1}^{II}$ to denote vertex set in which all the vertices in tree $\mathcal{T}$ are. And then, it is natural to group all newly added vertices by implementing $EUGM$ on every edge $uv$ in $\mathcal{T}$ into vertex set $\Lambda_{\mathcal{E}_{uv}}^{II}$ where symbol $\mathcal{E}_{uv}$ represents a specific path in stochastic uniform growth tree $\mathcal{T}_{II}(1)$ whose two end-vertices are adjacent in seed $\mathcal{T}$. The way to do this is to distinguish path $\mathcal{P}_{uv}$. These such paths $\mathcal{E}_{uv}$ are collected into a set $\mathcal{E}_{uv}^{II}$. Additionally, each vertex in $\Lambda_{\mathcal{E}_{uv}}^{II}$ is assigned a unique label $w_{uv}^{\alpha}$. By using these notations above, vertex set $V_{\mathcal{T}_{II}(1)}$ is expressed as $\bigcup_{uv\in\mathcal{T}}\Lambda_{\mathcal{E}_{uv}}^{II}\bigcup \Omega_{1}^{II}$. We are now ready to provide a rigorous proof to Eq.(\ref{Section-4-2-0}). As will be explained below, our computations are carried out in stages.

\emph{Case 1} Using $EUGM$, there will be $m_{i}$ vertices inserted into each edge $uv$ in seed $\mathcal{T}$ with probability $p_{i}$. As a consequence, it is not hard to see

\begin{equation}\label{Section-4-2-0-1}
\mathcal{W}_{\mathcal{T}_{II}(1)}(1)\triangleq\frac{1}{2}\sum_{u\in\Omega_{1}^{II}}\sum_{v\in\Omega_{1}^{II}}d'_{uv}=\mathcal{W}_{\mathcal{T}}\left(\sum_{i=1}^{\mu}p_{i}m_{i}+1\right),
\end{equation}
where summation $\sum_{i=1}^{\mu}p_{i}m_{i}+1$ will be replaced with symbol $\Phi_{4}^{II}$ in executing further arithmetics for the purpose of simplifying calculation.

\emph{Case 2} Following the aforementioned case, each of $m_{i}$ newly inserted vertices is attached $\nu$ paths having $n_{j}$ vertices with probability $q_{j}$ each. By utilizing a similar computational manner to that used to analyze case 2 of theorem 5, we can write

\begin{equation}\label{Section-4-2-0-2}
\begin{aligned}\mathcal{W}_{\mathcal{T}_{II}(1)}(2)&\triangleq\frac{1}{2}\sum_{\mathcal{E}_{uv}^{II}}\sum_{w_{uv}^{\alpha}\in\Lambda^{II}_{\mathcal{E}_{uv}}}\sum_{w_{uv}^{\beta}\in\Lambda^{II}_{\mathcal{E}_{uv}}}d'_{w_{uv}^{\alpha}w_{uv}^{\beta}}\\
&=(h-1)\sum_{i=1}^{\mu}p_{i}\sum_{l=2}^{m_{i}}\left(
                                                             \begin{array}{c}
                                                               l \\
                                                               2 \\
                                                             \end{array}
                                                           \right)\left(1+\nu\sum_{j=1}^{\nu}q_{j}n_{j}\right)^{2}+(h-1)\sum_{i=1}^{\mu}p_{i}m_{i}^{2}\left[\nu\sum_{j=1}^{\nu}q_{j}\left(
                                                             \begin{array}{c}
                                                               n_{j}+1 \\
                                                               2 \\
                                                             \end{array}
                                                           \right)\right]\\
                                                           &\quad+(h-1)\left[\sum_{i=1}^{\mu}p_{i}m_{i}\left(
                                                             \begin{array}{c}
                                                               \nu \\
                                                               2 \\
                                                             \end{array}
                                                           \right)+\nu^{2}\sum_{i=1}^{\mu}p_{i}\left(
                                                             \begin{array}{c}
                                                               m_{i} \\
                                                               2 \\
                                                             \end{array}
                                                           \right)\right]\left[2\sum_{j=1}^{\nu}q_{j}n_{j}\sum_{j=1}^{\nu}q_{j}\left(
                                                             \begin{array}{c}
                                                               n_{j}+1 \\
                                                               2 \\
                                                             \end{array}
                                                           \right)\right]\\
                                                           &\quad+(h-1)\sum_{i=1}^{\mu}p_{i}m_{i}\left[\nu\sum_{j=1}^{\nu}q_{j}\sum_{s=2}^{n_{j}}\left(
                                                             \begin{array}{c}
                                                               s \\
                                                               2 \\
                                                             \end{array}
                                                           \right)\right]
\end{aligned}.
\end{equation}
As before, the summation over all ``coefficients" of term $(h-1)$ in Eq.(\ref{Section-4-2-0-2}) is denoted by symbol $\Phi_{2}^{II}$ when we perform further computations.

\emph{Case 3} Now, we discuss the distance between vertex $u$ in set $\Omega_{1}^{II}$ and vertex $w_{xy}^{\alpha}$ in set $\bigcup_{xy\in\mathcal{T}}\Lambda_{\mathcal{E}_{xy}}^{II}$. Without loss of generality, we make use of $\mathcal{P}_{uw_{xy}^{\alpha}}$ to indicate that path linking vertex $u$ with $w_{xy}^{\alpha}$ in stochastic uniform growth tree $\mathcal{T}_{II}(1)$. In addition, path $\mathcal{P}_{uw_{xy}^{\alpha}}$ is assumed to possess two sub-paths $\mathcal{P}_{ux}$ and $\mathcal{P}_{xw_{xy}^{\alpha}}$. Note that we have supposed that vertex $y$ is always far away from vertex $u$ than vertex $x$. In this way, there is no influence on the future derivations. In what follows, we can encounter two sub-cases: (1) $w_{xy}^{\alpha}$ is some vertex inserted into edge $xy$ in seed $\mathcal{T}$, and (2) vertex $w_{xy}^{\alpha}$ is in some one of $\nu$ paths attached to some vertex $\theta$\footnote[3]{Here, we take symbol $\theta$ to indicate some vertex in set $\bigcup_{xy\in\mathcal{T}}\Lambda_{\mathcal{E}_{xy}}^{II}$, which is inserted into edge $xy$ in seed $\mathcal{T}$ through $EUGM$ directly, for convenience.} that is inserted into edge $xy$ in seed $\mathcal{T}$. In any subcase, we would like to replace sub-path $\mathcal{P}_{ux}$ with $\mathcal{P}_{uy}$ so as to build up a connection of quantity $\mathcal{W}_{\mathcal{T}_{II}(1)}(1)$ to $\mathcal{W}_{\mathcal{T}_{II}(1)}(3)$ as follows. For the first subcase, the distance $d'_{uw_{xy}^{\alpha}}$ satisfies $d'_{uw_{xy}^{\alpha}}=d'_{uy}-d'_{yw_{xy}^{\alpha}}.$
And, in the other subcase, the distance $d'_{uw_{xy}^{\alpha}}$ is given by $d'_{uw_{xy}^{\alpha}}=d'_{uy}-d'_{y\theta}+d'_{\theta w_{xy}^{\alpha}}.$
Based on the analysis above, we can obtain

\begin{equation}\label{Section-4-2-0-3}
\begin{aligned}\mathcal{W}_{\mathcal{T}_{II}(1)}(3)&\triangleq\sum_{u\in\Omega_{1}^{II}}\sum_{\mathcal{E}_{xy}^{II}}\sum_{w_{xy}^{\alpha}\in\Lambda^{II}_{xy}}d'_{uw_{xy}^{\alpha}}\\
&=2\mathcal{W}_{\mathcal{T}_{II}(1)}(1)\left(\sum_{i=1}^{\mu}p_{i}m_{i}\right)\left(\nu\sum_{j=1}^{\nu}q_{j}n_{j}+1\right)\\
&\quad-2\left(
                                                             \begin{array}{c}
                                                              h \\
                                                               2 \\
                                                             \end{array}
                                                           \right)\left(\sum_{i=1}^{\mu}p_{i}m_{i}+1\right)\left(\sum_{i=1}^{\mu}p_{i}m_{i}\right)\left(\nu\sum_{j=1}^{\nu}q_{j}n_{j}+1\right)\\
&\quad+2\left(
                                                             \begin{array}{c}
                                                              h \\
                                                               2 \\
                                                             \end{array}
                                                           \right)\left[\sum_{i=1}^{\mu}p_{i}\left(
                                                             \begin{array}{c}
                                                               m_{i}+1 \\
                                                               2 \\
                                                             \end{array}
                                                           \right)\left(\nu\sum_{j=1}^{\nu}q_{j}n_{j}+1\right)+\sum_{i=1}^{\mu}p_{i}m_{i}\nu\sum_{j=1}^{\nu}q_{i}\left(
                                                             \begin{array}{c}
                                                               n_{j}+1 \\
                                                               2 \\
                                                             \end{array}
                                                           \right)\right]
\end{aligned}.
\end{equation}
For ease of exposition, we still need to introduce two symbols, $\Phi_{1}^{II}$ and $\Phi_{3}^{II}$, as stated early. More specifically,

$$\Phi_{1}^{II}=\left(\sum_{i=1}^{\mu}p_{i}m_{i}\right)\left(\nu\sum_{j=1}^{\nu}q_{j}n_{j}+1\right),$$
along with
$$\Phi_{3}^{II}=\sum_{i=1}^{\mu}p_{i}\left(
                                                             \begin{array}{c}
                                                               m_{i}+1 \\
                                                               2 \\
                                                             \end{array}
                                                           \right)\left(\nu\sum_{j=1}^{\nu}q_{j}n_{j}+1\right)+\sum_{i=1}^{\mu}p_{i}m_{i}\nu\sum_{j=1}^{\nu}q_{i}\left(
                                                             \begin{array}{c}
                                                               n_{j}+1 \\
                                                               2 \\
                                                             \end{array}
                                                           \right).$$

\emph{Case 4} The last task is to derive analytic solution to summation $\mathcal{W}_{\mathcal{T}_{II}(1)}(4)$ over distances between vertices $w_{uv}^{\alpha}$ and $w_{xy}^{\beta}$ in tree $\mathcal{T}_{II}(1)$ where $uv$ is not identical to $xy$. Along the similar thought in cases 2 and 3, we omit detailed demonstration and straightforwardly provide

\begin{equation}\label{Section-4-2-0-4}
\begin{aligned}\mathcal{W}_{\mathcal{T}_{II}(1)}(4)&\triangleq\sum_{uv(\neq xy)\in\mathcal{T}}\sum_{w_{uv}^{\alpha}\in\Lambda^{II}_{xy}}\sum_{w_{xy}^{\beta}\in\Lambda^{II}_{xy}}d'_{w_{uv}^{\alpha}w_{xy}^{\beta}}\\
&=\sum_{uv(\neq xy)\in\mathcal{T}}\sum_{w_{uv}^{\alpha}\in\Lambda^{II}_{xy}}\sum_{w_{xy}^{\beta}\in\Lambda^{II}_{xy}}\left[d'_{uy}-2\left(\sum_{i=1}^{\mu}p_{i}m_{i}+1\right)+d'_{vw_{uv}^{\alpha}}+d'_{xw_{xy}^{\beta}}\right]\\
&=\left[\mathcal{W}_{\mathcal{T}_{II}(1)}(1)-(h-1)\Phi_{4}^{II}\right]\left(\Phi_{1}^{II}\right)^{2}-2\Phi_{4}^{II}\left(\Phi_{1}^{II}\right)^{2}\left(
                                                             \begin{array}{c}
                                                               h-1 \\
                                                               2 \\
                                                             \end{array}
                                                           \right)+2\Phi_{4}^{II}\Phi_{3}^{II}\left(
                                                             \begin{array}{c}
                                                              h-1 \\
                                                               2 \\
                                                             \end{array}
                                                           \right)
\end{aligned},
\end{equation}
in which we have used some symbols that have the same meanings as above.

So far, substituting the results from Eqs.(\ref{Section-4-2-0-1})-(\ref{Section-4-2-0-4}) into expression $\mathcal{W}_{\mathcal{T}_{II}(1)}=\sum_{i=1}^{4}\mathcal{W}_{\mathcal{T}_{II}(1)}(i)$ yields the same consequence as in Eq.(\ref{Section-4-2-0}), implying that theorem 9 is complete. $\hfill\Box$

Upon an arbitrary tree $\mathcal{T}$ as seed, the final graph $\mathcal{T}_{II}(t)$ is recursively constructed via executing $EUGM$
$t$ steps. After that, the vertex number $|\mathcal{V}_{\mathcal{T}_{II}(t)}|$ is easy to calculate in an iterative way, as below

\begin{equation}\label{Section-4-2-1}
|\mathcal{V}_{\mathcal{T}_{II}(t)}|=(h-1)\left(1+\sum_{i=1}^{\mu}p_{i}m_{i}+\sum_{i=1}^{\mu}p_{i}m_{i}\nu\sum_{j=1}^{\nu}q_{j}n_{j}\right)^{t}+1.
\end{equation}

From now on, let us focus on the calculation of Wiener index $\mathcal{W}_{\mathcal{T}_{II}(t)}$ of stochastic uniform growth tree $\mathcal{T}_{II}(t)$. As stated previously, this issue can also be effortlessly addressed by calculating recurrence relation in terms of Eqs.(\ref{Section-4-2-0}) and (\ref{Section-4-2-1}). Thus, we omit the detailed derivation and attach the final formula in the following proposition.

\textbf{Proposition 10} The solution of Wiener index $\mathcal{W}_{\mathcal{T}_{II}(t)}$ of tree $\mathcal{T}_{II}(t)$ is given by

\begin{equation}\label{Section-4-2-2}
\begin{aligned}\mathcal{W}_{\mathcal{T}_{II}(t)}&=\mathcal{W}_{\mathcal{T}}\left[\Phi_{4}^{II}\left(1+\Phi_{1}^{II}\right)^{2}\right]^{t}+\left(\Psi_{1}^{III}+\Psi_{2}^{III}+\Psi_{3}^{III}\right)\frac{\left[\Phi_{4}^{II}\left(1+\Phi_{1}^{II}\right)^{2}\right]^{t}-1}{\Phi_{4}^{II}\left(1+\Phi_{1}^{II}\right)^{2}-1}\\
&+(h-1)\left(2\Psi_{1}^{II}+\Psi_{2}^{II}\right)\left(1+\Phi_{1}^{II}\right)^{t-1}\frac{\left[\Phi_{4}^{II}\left(1+\Phi_{1}^{II}\right)\right]^{t}-1}{\Phi_{4}^{II}\left(1+\Phi_{1}^{II}\right)-1}\\
&+\Psi_{1}^{II}(h-1)^{2}\left(1+\Phi_{1}^{II}\right)^{2(t-1)}\frac{\left(\Phi_{4}^{II}\right)^{t}-1}{\Phi_{4}^{II}-1}\\
\end{aligned},
\end{equation}
in which we have taken advantage of three additional symbols $\Psi_{i}^{II}$ for convenience. Their own specific implications are as follows

$$\Psi_{1}^{II}=\Phi_{1}^{II}\Phi_{3}^{II}-\Phi_{4}^{II}\left(\Phi_{1}^{II}\right)^{2}-\Phi_{4}^{II}\Phi_{1}^{II}+\Phi_{3}^{II},$$

$$\Psi_{2}^{II}=2\Phi_{4}^{II}\left(\Phi_{1}^{II}\right)^{2}-3\Phi_{1}^{II}\Phi_{3}^{II}+\Phi_{4}^{II}\Phi_{1}^{II}+\Phi_{2}^{II}-\Phi_{3}^{II},\quad \text{and,}\quad \Psi_{3}^{II}=2\Phi_{1}^{II}\Phi_{3}^{II}-\Phi_{4}^{II}\left(\Phi_{1}^{II}\right)^{2}-\Phi_{2}^{II}.$$
Here, there is a little surprise, i.e.,
\begin{equation}\label{Section-4-2-2-1}
\Psi_{1}^{II}+\Psi_{2}^{II}+\Psi_{3}^{II}=0.
\end{equation}
That is to say, we only need to derive arbitrary two of parameters $\Psi_{i}^{II}$  when calculating the $\mathcal{W}$-polynomial for tree $\mathcal{W}_{\mathcal{T}_{II}(t)}$ with respect of parameters $\mathcal{W}_{\mathcal{T}}$ and $h$ on seed $\mathcal{T}$.

As a special member in stochastic uniform growth tree family $\mathcal{T}_{II}(1)$, the $m$-th order subdivision tree $\mathcal{T}^{m}$ is deterministic, and may be conveniently produced by setting parameters $m_{i}=m$ for all $i\in[1,\mu]$ and $\nu=0$ in $EUGM$. This kind of trees have been commonly-studied in graph theory \cite{Bondy-2008}. Here, for our purpose, the closed-form solution to corresponding Wiener index on trees of such type can be immediately obtained from Eq.(\ref{Section-4-2-0}).

\textbf{Corollary 11} The closed-form solution to Wiener index $\mathcal{W}_{\mathcal{T}^{m}}$ on $m$-th order subdivision tree $\mathcal{T}^{m}$ is given in the following form

\begin{equation}\label{Section-4-2-3}
\mathcal{W}_{\mathcal{T}^{m}}=(m+1)^{3}\mathcal{W}_{\mathcal{T}}-\frac{m(m+1)^{2}}{2}h^{2}+\left[\frac{m(m+1)^{2}}{2}+\sum_{l=2}^{m}\left(
                                                             \begin{array}{c}
                                                               l \\
                                                               2 \\
                                                             \end{array}
                                                           \right)\right]h-\sum_{l=2}^{m}\left(
                                                             \begin{array}{c}
                                                               l \\
                                                               2 \\
                                                             \end{array}
                                                           \right).
\end{equation}

In particular, when setting $m=1$, the $m$-th order subdivision tree $\mathcal{T}^{m}$ is reduced as the subdivision tree $\mathcal{T}'$. Then, the corresponding Wiener index $\mathcal{W}_{\mathcal{T}'}$ is calculated to equal

$$\mathcal{W}_{\mathcal{T}'}=8\mathcal{W}_{\mathcal{T}}-2h(h-1),$$
which is identical to result in our previous work \cite{Ma-2020} (see theorem 1 in \cite{Ma-2020} for more details).

We are now in a position where the formula of mean first-passage time $\overline{\mathcal{F}}_{\mathcal{T}_{II}(t)}$ on stochastic uniform growth tree $\mathcal{T}_{II}(t)$ need be derived analytically. In practice, this can be derived using results in Eqs.(\ref{Section-4-2-1})  and (\ref{Section-4-2-2}) by virtue of statement from Eq.(\ref{Section-3-4-0}). For brevity and convenience, we omit the concrete derivation, and the final expression is written as below.

\textbf{Theorem 12} The analytic solution to mean first-passage time $\overline{\mathcal{F}}_{\mathcal{T}_{II}(t)}$ for random walks on tree $\mathcal{T}_{II}(t)$ is given by

\begin{equation}\label{Section-4-2-4}
\begin{aligned}\overline{\mathcal{F}}_{\mathcal{T}_{II}(t)}&=\mathcal{W}_{\mathcal{T}}\frac{2\left[\Phi_{4}^{II}\left(1+\Phi_{1}^{II}\right)^{2}\right]^{t}}{(h-1)\left(1+\Phi_{1}^{II}\right)^{t}+1}+\frac{2\left(\Psi_{1}^{III}+\Psi_{2}^{III}+\Psi_{3}^{III}\right)}{(h-1)\left(1+\Phi_{1}^{II}\right)^{t}+1}\times\frac{\left[\Phi_{4}^{II}\left(1+\Phi_{1}^{II}\right)^{2}\right]^{t}-1}{\Phi_{4}^{II}\left(1+\Phi_{1}^{II}\right)^{2}-1}\\
&+2(h-1)\frac{\left(2\Psi_{1}^{II}+\Psi_{2}^{II}\right)\left(1+\Phi_{1}^{II}\right)^{t-1}}{(h-1)\left(1+\Phi_{1}^{II}\right)^{t}+1}\times\frac{\left[\Phi_{4}^{II}\left(1+\Phi_{1}^{II}\right)\right]^{t}-1}{\Phi_{4}^{II}\left(1+\Phi_{1}^{II}\right)-1}\\
&+2(h-1)^{2}\frac{\Psi_{1}^{II}\left(1+\Phi_{1}^{II}\right)^{2(t-1)}}{(h-1)\left(1+\Phi_{1}^{II}\right)^{t}+1}\times\frac{\left(\Phi_{4}^{II}\right)^{t}-1}{\Phi_{4}^{II}-1}\\
\end{aligned}.
\end{equation}

Besides that, in the limit of large graph size, the result above will behave a power scaling over variable $|V_{\mathcal{T}_{II}(t)}|$,

$$\overline{\mathcal{F}}_{\mathcal{T}_{II}(t)}=O\left(|V_{\mathcal{T}_{II}(t)}|^{\Theta_{II}}\right),\qquad \Theta_{II}=1+\frac{\ln\left(\sum_{i=1}^{\mu}p_{i}m_{i}+1\right)}{\ln\left[1+\left(\sum_{i=1}^{\mu}p_{i}m_{i}\right)\left(\nu\sum_{j=1}^{\nu}q_{j}n_{j}+1\right)\right] }.$$
Obviously, parameter $\Theta_{II}$ is no larger than $2$. Yet, it asymptotically tends constant $2$ when assuming $\nu=0$.

It has been shown in Section 3 that the famous $T$-graph $\mathcal{Y}_{II}(t)$ is the simplest member in stochastic uniform growth tree $\mathcal{T}_{II}(t)$. Some intriguing structural parameters planted on $\mathcal{Y}_{II}(t)$ including mean first-passage time $\overline{\mathcal{F}}_{\mathcal{Y}_{II}(t)}$ have been studied using other methods in the past \cite{Agliari-2008} (see Eq.(13) in \cite{Agliari-2008} for more details). Here, we only need to substitute some initial conditions, namely, $\mathcal{W}_{\mathcal{T}}=1$, $h=2$, $m_{i}=1$ for $i\in[1,\mu]$, $\nu=1$ as well as $n_{1}=1$, into Eq.(\ref{Section-4-2-4}) in order to obtain the corresponding formula for quantity $\overline{\mathcal{F}}_{\mathcal{Y}_{II}(t)}$.

\textbf{Corollary 13} The exact solution to mean first-passage time $\overline{\mathcal{F}}_{\mathcal{Y}_{II}(t)}$ on the well-known $T$-graph $\mathcal{Y}_{II}(t)$ is

\begin{equation}\label{Section-4-2-5}
\overline{\mathcal{F}}_{\mathcal{Y}_{II}(t)}=\frac{2}{3^{t}+1}\left(18^{t}-2\times\frac{18^{t}-3^{t}}{5}-\frac{18^{t}-9^{t}}{3}\right).
\end{equation}

To make further progress, if we suppose that in stochastic uniform growth tree $\mathcal{T}_{II}(t)$, the seed is still an edge and parameters $m_{i}=1$, $n_{j}=1$ for all $j$ ($j\in[1,\nu]$), then the resulting deterministic graph is $\nu$-fractal tree $\mathcal{T}_{\nu}(t)$. In \cite{Lin-2010}, the mean first-passage time $\overline{\mathcal{F}}_{\mathcal{T}_{\nu}(t)}$ on tree $\mathcal{T}_{\nu}(t)$ has been reported using spectral method (see Eq.(64) in \cite{Lin-2010} for more details). On the other hand, the corresponding formula $\overline{\mathcal{F}}_{\mathcal{T}_{\nu}(t)}$ is able to be exactly obtained through substituting parameters which are related to tree $\mathcal{T}_{\nu}(t)$ into Eq.(\ref{Section-4-2-4}), which is as follows

$$\overline{\mathcal{F}}_{\mathcal{T}_{\nu}(t)}=\frac{2}{(\nu+2)^{t}+1}\left\{\left[2(\nu+1)^{2}\right]^{t}-(\nu+1)(\nu+2)^{t}\frac{(2\nu+4)^{t}-1}{2\nu+3}-(\nu+2)^{2t-1}\left(2^{t}-1\right)\right\}.$$

\subsection{Tree $\mathcal{T}_{III}(t)$}

\textbf{Theorem 14} The solution of Wiener index $\mathcal{W}_{\mathcal{T}_{III}(1)}$ of tree $\mathcal{T}_{III}(1)$ is

\begin{equation}\label{Section-4-3-0}
\begin{aligned}\mathcal{W}_{\mathcal{T}_{III}(1)}
&=\mathcal{W}_{\mathcal{T}}\left(2\sum_{i=1}^{\mu}p_{i}m_{i}+1\right)\left(\mu\sum_{i=1}^{\mu}p_{i}m_{i}+1\right)^{2}+(\mu-2)\left(\mu\sum_{i=1}^{\mu}p_{i}m_{i}+1\right)\sum_{i=1}^{\mu}p_{i}\left(
                                                             \begin{array}{c}
                                                               m_{i}+1 \\
                                                               2 \\
                                                             \end{array}
                                                           \right)h^{2}\\
&+\left\{\left(\mu\sum_{i=1}^{\mu}p_{i}m_{i}+2\right)\sum_{i}^{\mu}p_{i}\left(
                                                             \begin{array}{c}
                                                               m_{i}+1 \\
                                                               2 \\
                                                             \end{array}
                                                           \right)+\mu\sum_{i}^{\mu}p_{i}\sum_{l=2}^{m_{i}}\left(
                                                             \begin{array}{c}
                                                               l \\
                                                               2 \\
                                                             \end{array}
                                                           \right)\right\}h
\end{aligned}.
\end{equation}
Similarly, we can write the $\mathcal{W}$-polynomial for tree $\mathcal{T}_{III}(1)$ in the following form

\begin{equation}\label{Section-4-3-0-0}
f_{III}(\mathcal{W}_{\mathcal{T}},h)\triangleq a_{III}\mathcal{W}_{\mathcal{T}}+b_{III}h^{2}+c_{II}h.
\end{equation}

\textbf{\emph{Proof}} In fact, there exist some similarities between $VUGM$ and $MUGM$. For instance, a star-like subgraph with $\mu$ ``tentacles" is to be created for each vertex $u$ in seed $\mathcal{T}$ when performing operation on vertex $u$. More specifically, this type of star-like subgraph contains vertex $u$ as center and those newly added vertices in terms of vertex $u$. The latter vertices are grouped into set $\Lambda_{u}^{III}$ for convenience and our purpose. As above, each vertex in set $\Lambda_{u}^{III}$ is remarked by a unique symbol $u_{\alpha}$. Additionally, we still make use of $\bigcup_{u\in\mathcal{T}}\Lambda_{u}^{III}$ to represent set composed of all the new vertices introduced into tree $\mathcal{T}_{III}(1)$ when applying $MUGM$ to every vertex in seed $\mathcal{T}$, and define notation $\Omega_{1}^{III}=\mathcal{V}_{\mathcal{T}_{III}(1)}/\bigcup_{u\in\mathcal{T}}\Lambda_{u}^{III}$ to be vertex set to which all the vertices in seed $\mathcal{T}$ belong. Now, let us divert more attention to calculation of Wiener index of tree $\mathcal{T}_{III}(1)$. This is dealt with using a similar method as recommended previously. At the same time, it is noteworthy that some formulae will be shown straightforwardly without detailed description according to the same derivation in the development of validating Eq.(\ref{Section-4-1-0}). Reader refers subsection 4.1 for more details.

\emph{Case 1} By definition, it is clear to see that there will be $2\sum_{i=1}^{\mu}p_{i}m_{i}$ vertices inserted into each edge $uv$ in seed $\mathcal{T}$. This implies that the sum $\mathcal{W}_{\mathcal{T}_{III}(1)}(1)$ over all distances between two arbitrarily distinct vertices in set $\Omega_{1}^{III}$ is certainly subject to the following formula

\begin{equation}\label{Section-4-3-0-1}
\mathcal{W}_{\mathcal{T}_{III}(1)}(1)\triangleq\frac{1}{2}\sum_{u\in\Omega_{1}^{III}}\sum_{v\in\Omega_{1}^{III}}d'_{uv}=\mathcal{W}_{\mathcal{T}}\left(2\sum_{i=1}^{\mu}p_{i}m_{i}+1\right).
\end{equation}

\emph{Case 2} For both vertex in set $\Lambda_{u}^{III}$ and vertex $u$ in set $\Omega_{1}^{III}$, we surely have

\begin{equation}\label{Section-4-3-0-2}
\begin{aligned}\mathcal{W}_{\mathcal{T}_{III}(1)}(2)&\triangleq\sum_{u\in\Omega_{1}^{III}}\left(\sum_{u_{\alpha }\in\Lambda_{u}^{III}}d'_{uu_{\alpha}}+\frac{1}{2}\sum_{u_{\alpha}\in\Lambda_{u}^{III}}\sum_{u_{\beta}\in\Lambda_{u}^{III}}d'_{u_{\alpha}u_{\beta}}\right)\\
&=h\left[\mu\sum_{i=1}^{\mu}p_{i}\left(
                                                             \begin{array}{c}
                                                               m_{i}+1 \\
                                                               2 \\
                                                             \end{array}
                                                           \right)+\sum_{i=1}^{\mu}2p_{i}m_{i}\left(
                                                             \begin{array}{c}
                                                               \mu \\
                                                               2 \\
                                                             \end{array}
                                                           \right)\sum_{i=1}^{\mu}p_{i}\left(
                                                             \begin{array}{c}
                                                               m_{i}+1 \\
                                                               2 \\
                                                             \end{array}
                                                           \right)+\mu\sum_{i}^{\mu}p_{i}\sum_{l=2}^{m_{i}}\left(
                                                             \begin{array}{c}
                                                               l \\
                                                               2 \\
                                                             \end{array}
                                                           \right)\right]
                                                          \end{aligned}.
\end{equation}
This is completely the same as Eq.(\ref{Section-4-1-0-2}).

\emph{Case 3} Apart from some similarities between $VUGM$ and $MUGM$, there are a few differences as will be stated shortly. Given a pair of vertices, say $u$ and $v_{\alpha}$ where $u$ differs for $v$, we denote by $\mathcal{P}_{uv_{\alpha}}$ path connecting vertex $u$ to $v_{\alpha}$ in stochastic uniform growth tree $\mathcal{T}_{III}(1)$. In view of $MUGM$, we confirm that each of $\mu-1$ ``tentacles" in star-like subgraph whose center is vertex $v$ can be referred to as an expansion of path $\mathcal{P}_{uv}$. To put this another way, path $\mathcal{P}_{uv_{\alpha}}$ is based on paths $\mathcal{P}_{uv}$ and $\mathcal{P}_{vv_{\alpha}}$ via conjunction on vertex $v$. This leads to a relation $d'_{uv_{\alpha}}=d'_{uv}+d'_{vv}.$
On the other hand, the left ``tentacle" in star-like subgraph is in fact a contraction of path $\mathcal{P}_{uv}$ itself. Specifically speaking, path $\mathcal{P}_{uv_{\alpha}}$ is obtained from  path $\mathcal{P}_{uv}$ by deleting $\mathcal{P}_{vv_{\alpha}}$, which results in the next expression $d'_{uv_{\alpha}}=d'_{uv}-d'_{vv}.$
Taken together, we derive the solution to summation $\mathcal{W}_{\mathcal{T}_{III}(1)}(3)$ over distances of all possible vertex pairs $u$ and $v_{\alpha}$ of this kind in tree $\mathcal{T}_{III}(1)$, as follows

\begin{equation}\label{Section-4-3-0-3}
\begin{aligned}\mathcal{W}_{\mathcal{T}_{III}(1)}(3)&\triangleq\frac{1}{2}\sum_{u\in\Omega_{1}^{III}}\sum_{v(\neq u)\in\Omega_{1}^{III}}\sum_{v_{\alpha}\in\Lambda_{v}^{III}}d'_{uv_{\alpha}}\\
&=2\mathcal{W}_{\mathcal{T}_{III}(1)}(1)\left(\mu\sum_{i=1}^{\mu}p_{i}m_{i}\right)+2(\mu-2)\left(
                                                             \begin{array}{c}
                                                               h \\
                                                               2 \\
                                                             \end{array}
                                                           \right)\sum_{i=1}^{\mu}p_{i}\left(
                                                             \begin{array}{c}
                                                               m_{i}+1 \\
                                                               2 \\
                                                             \end{array}
                                                           \right)
\end{aligned}.
\end{equation}

\emph{Case 4} Finally, let us evaluate the contribution from distance between arbitrary pair of vertices $u_{\alpha}$ and $v_{\beta}$ to Wiener index $\mathcal{W}_{\mathcal{T}_{III}(1)}$, which is thought of as $\mathcal{W}_{\mathcal{T}_{III}(1)}(4)$, i.e.,

$$\mathcal{W}_{\mathcal{T}_{III}(1)}(4)\triangleq\sum_{u\in\Omega_{1}^{III}}\sum_{v(\neq u)\in\Omega_{1}^{III}}\sum_{u_{\alpha}\in\Lambda_{u}^{III}}\sum_{v_{\beta}\in\Lambda_{v}^{III}}d'_{u_{\alpha}v_{\beta}.}$$
By analogy with demonstration in previous cases, we can determine the analytic solution to quantity $\mathcal{W}_{\mathcal{T}_{III}(1)}(4)$ by first considering two corresponding central vertices $u$ and $v$ with respect to a given pair of vertices $u_{\alpha}$ and $v_{\beta}$. Specifically, for star-like subgraph whose center is $u$, there must be $\mu-1$ ``tentacles" as expansions of path $\mathcal{P}_{uv}$, each being towards the outside along the direction from end-vertex $v$ to $u$, and $1$ ``tentacle" as contraction of path $\mathcal{P}_{uv}$ towards the inside along the opposite direction, namely, orientation from end-vertex $u$ to $v$. Taking into account nature of $MUGM$, there are four different combinatorial manners in situation mentioned above by means of both expansion and contraction along two distinct directions. For the sake of argument, we omit concrete calculation for each combinatorial manner and immediately write

\begin{equation}\label{Section-4-3-0-4}
\begin{aligned}\mathcal{W}_{\mathcal{T}_{III}(1)}(4)&=\mathcal{W}_{\mathcal{T}_{III}(1)}(1)\left((\mu-1)\sum_{i=1}^{\mu}p_{i}m_{i}\right)^{2}+2(\mu-1)\sum_{i=1}^{\mu}p_{i}m_{i}\left(
                                                             \begin{array}{c}
                                                               h \\
                                                               2 \\
                                                             \end{array}
                                                           \right)\left[(\mu-1)\sum_{i=1}^{\mu}p_{i}\left(
                                                             \begin{array}{c}
                                                               m_{i}+1 \\
                                                               2 \\
                                                             \end{array}
                                                           \right)\right]\\
                                                           &\quad+2\mathcal{W}_{\mathcal{T}_{III}(1)}(1)(\mu-1)\left(\sum_{i=1}^{\mu}p_{i}m_{i}\right)^{2}-2\left(
                                                             \begin{array}{c}
                                                               h \\
                                                               2 \\
                                                             \end{array}
                                                           \right)(\mu-1)\sum_{i=1}^{\mu}p_{i}m_{i}\sum_{i=1}^{\mu}p_{i}\left(
                                                             \begin{array}{c}
                                                               m_{i}+1 \\
                                                               2 \\
                                                             \end{array}
                                                           \right)\\
                                                           &\quad+2\sum_{i=1}^{\mu}p_{i}m_{i}\left(
                                                             \begin{array}{c}
                                                               h \\
                                                               2 \\
                                                             \end{array}
                                                           \right)\left[(\mu-1)\sum_{i=1}^{\mu}p_{i}\left(
                                                             \begin{array}{c}
                                                               m_{i}+1 \\
                                                               2 \\
                                                             \end{array}
                                                           \right)\right]\\
                                                           &\quad+\mathcal{W}_{\mathcal{T}_{III}(1)}(1)\left(\sum_{i=1}^{\mu}p_{i}m_{i}\right)^{2}-2\left(
                                                             \begin{array}{c}
                                                               h \\
                                                               2 \\
                                                             \end{array}
                                                           \right)\left[\left(\sum_{i=1}^{\mu}p_{i}m_{i}\right)\sum_{i=1}^{\mu}p_{i}\left(
                                                             \begin{array}{c}
                                                               m_{i}+1 \\
                                                               2 \\
                                                             \end{array}
                                                           \right)\right]
\end{aligned}.
\end{equation}
Using some fundamental arithmetics, the above equation simplifies to

\begin{equation}\label{Section-4-3-0-5}
\mathcal{W}_{\mathcal{T}_{III}(1)}(4)=\mathcal{W}_{\mathcal{T}_{III}(1)}(1)\left(\mu\sum_{i=1}^{\mu}p_{i}m_{i}\right)^{2}+2(\mu^{2}-2\mu)\sum_{i=1}^{\mu}p_{i}m_{i}\left(
                                                             \begin{array}{c}
                                                               h \\
                                                               2 \\
                                                             \end{array}
                                                           \right)\left[\sum_{i=1}^{\mu}p_{i}\left(
                                                             \begin{array}{c}
                                                               m_{i}+1 \\
                                                               2 \\
                                                             \end{array}
                                                           \right)\right].
\end{equation}

Putting all things together yields the same result as in Eq.(\ref{Section-4-3-0}), which suggests that we complete the proof to theorem 14. $\hfill\Box$

From Eq.(\ref{Section-4-3-0}), we can find that there is in essence a recursive relation of both Wiener index $\mathcal{W}_{\mathcal{T}_{III}(1)}$ and $\mathcal{W}_{\mathcal{T}}$. This means that the analytic solution to Wiener index $\mathcal{W}_{\mathcal{T}_{III}(t)}$ may be easily obtained in an iterative fashion after estimating the total number of vertices in stochastic uniform growth tree $\mathcal{T}_{III}(t)$. To this end, in view of the specific growth way to construct tree $\mathcal{T}_{III}(t)$, the vertex number $|\mathcal{V}_{\mathcal{T}_{III}(t)}|$ follows

\begin{equation}\label{Section-4-3-1}
|\mathcal{V}_{\mathcal{T}_{III}(t)}|=h\left(1+\mu\sum_{i=1}^{\mu}p_{i}m_{i}\right)^{t}.
\end{equation}
Clearly, this is identical to the vertex number of tree $\mathcal{T}_{I}(t)$. This is another similarity between both types of trees,
$\mathcal{T}_{I}(t)$ and $\mathcal{T}_{III}(t)$. As before, an iterative calculation firmly produces the analytic formula of $\mathcal{W}_{\mathcal{T}_{III}(t)}$.

\textbf{Proposition 15} The solution of Wiener index $\mathcal{W}_{\mathcal{T}_{III}(t)}$ of tree $\mathcal{T}_{III}(t)$ is

\begin{equation}\label{Section-4-3-2}
\begin{aligned}\mathcal{W}_{\mathcal{T}_{III}(t)}
&=\mathcal{W}_{\mathcal{T}}\left(2\sum_{i=1}^{\mu}p_{i}m_{i}+1\right)^{t}\left(\mu\sum_{i=1}^{\mu}p_{i}m_{i}+1\right)^{2t}+(\mu-2)h^{2}\sum_{i=1}^{\mu}p_{i}\left(
                                                             \begin{array}{c}
                                                               m_{i}+1 \\
                                                               2 \\
                                                             \end{array}
                                                           \right)\left(\mu\sum_{i=1}^{\mu}p_{i}m_{i}+1\right)^{2t-1}\Phi^{III}_{1}\\
&+h\left\{\left(\mu\sum_{i=1}^{\mu}p_{i}m_{i}+2\right)\sum_{i}^{\mu}p_{i}\left(
                                                             \begin{array}{c}
                                                               m_{i}+1 \\
                                                               2 \\
                                                             \end{array}
                                                           \right)+\mu\sum_{i}^{\mu}p_{i}\sum_{l=2}^{m_{i}}\left(
                                                             \begin{array}{c}
                                                               l \\
                                                               2 \\
                                                             \end{array}
                                                           \right)\right\}\left(\mu\sum_{i=1}^{\mu}p_{i}m_{i}+1\right)^{t-1}\Phi^{III}_{2}
\end{aligned},
\end{equation}
where we make use of symbols $\Phi^{III}_{1}$ and $\Phi^{III}_{2}$, which are given in the following form

$$\Phi^{III}_{1}=\frac{\left(2\sum_{i=1}^{\mu}p_{i}m_{i}+1\right)^{t}-1}{2\sum_{i=1}^{\mu}p_{i}m_{i}},\quad \Phi^{III}_{2}=\frac{\left[\left(\mu\sum_{i=1}^{\mu}p_{i}m_{i}+1\right)\left(2\sum_{i=1}^{\mu}p_{i}m_{i}+1\right)\right]^{t}-1}{\left(\mu\sum_{i=1}^{\mu}p_{i}m_{i}+1\right)\left(2\sum_{i=1}^{\mu}p_{i}m_{i}+1\right)-1}.$$

For brevity and convenience, we omit more details about derivation. At meantime, based on Eq.(\ref{Section-3-4-0}), we convert the result in Eq.(\ref{Section-4-3-2}) into expression of mean first-passage time on stochastic uniform growth tree $\mathcal{T}_{III}(t)$, which is as below.

\textbf{Theorem 16} The analytic solution to mean first-passage time $\overline{\mathcal{F}}_{\mathcal{T}_{III}(t)}$ for random walks on tree $\mathcal{T}_{III}(t)$ is given by

\begin{equation}\label{Section-4-3-3}
\begin{aligned}\overline{\mathcal{F}}_{\mathcal{T}_{III}(t)}
&=\frac{2\mathcal{W}_{\mathcal{T}}}{h}\left(2\sum_{i=1}^{\mu}p_{i}m_{i}+1\right)^{t}\left(\mu\sum_{i=1}^{\mu}p_{i}m_{i}+1\right)^{t}+2(\mu-2)h\sum_{i=1}^{\mu}p_{i}\left(
                                                             \begin{array}{c}
                                                               m_{i}+1 \\
                                                               2 \\
                                                             \end{array}
                                                           \right)\left(\mu\sum_{i=1}^{\mu}p_{i}m_{i}+1\right)^{t-1}\Phi^{III}_{1}\\
&+2\left\{\left(\mu\sum_{i=1}^{\mu}p_{i}m_{i}+2\right)\sum_{i}^{\mu}p_{i}\left(
                                                             \begin{array}{c}
                                                               m_{i}+1 \\
                                                               2 \\
                                                             \end{array}
                                                           \right)+\mu\sum_{i}^{\mu}p_{i}\sum_{l=2}^{m_{i}}\left(
                                                             \begin{array}{c}
                                                               l \\
                                                               2 \\
                                                             \end{array}
                                                           \right)\right\}\left(\mu\sum_{i=1}^{\mu}p_{i}m_{i}+1\right)^{-1}\Phi^{III}_{2}
\end{aligned}.
\end{equation}

More generally, we are interested in the scaling behavior of quantity $\overline{\mathcal{F}}_{\mathcal{T}_{III}(t)}$ in the large graph size limit. On the basis of Eq.(\ref{Section-4-3-1}), when considering case $t\rightarrow \infty$, there is a relationship

$$\overline{\mathcal{F}}_{\mathcal{T}_{III}(t)}=O\left(|\mathcal{V}_{\mathcal{T}_{III}(t)}|^{\Theta_{III}}\right), \quad \Theta_{III}=1+\frac{\ln \left(2\sum_{i=1}^{\mu}p_{i}m_{i}+1\right)}{\ln\left(\mu\sum_{i=1}^{\mu}p_{i}m_{i}+1\right)},$$
suggesting that $\overline{\mathcal{F}}_{\mathcal{T}_{III}(t)}$ is not linearly correlated with vertex number $|\mathcal{V}_{\mathcal{T}_{III}(t)}|$ but a power function as parameter $|\mathcal{V}_{\mathcal{T}_{III}(t)}|$ in form. Meanwhile, the power exponent $\Theta_{III}$ is strictly less than constant $2$ and approaches unit gradually as $\mu\rightarrow \infty$.

As a case study, let us revisit the famous Vicsek fractal $V_{1}^{\mu}(t)$ in which the seed $\mathcal{T}$ is a star $S_{\mu}$, namely, $\mathcal{W}_{\mathcal{T}}=\mu^{2}$ and $h=|\mathcal{V}_{\mathcal{T}}|=\mu+1$. Substituting these pre-designated parameters into Eq.(\ref{Section-4-3-3}) yields the previously published result in \cite{Zhang-2010} (see Eq.(20) in \cite{Zhang-2010} for more details). To make our statement more self-contained, the closed-form solution is written in the next corollary.

\textbf{Corollary 17} The exact solution to mean first-passage time $\overline{\mathcal{F}}_{V_{1}^{\mu}(t)}$ on the famous Vicsek fractal $V_{1}^{\mu}(t)$ is

\begin{equation}\label{Section-4-3-4}
\overline{\mathcal{F}}_{V_{1}^{\mu}(t)}=\frac{2\mu^{2}}{\mu+1}(3\mu+3)^{t}+(\mu-2)(\mu+1)^{t}(3^{t}-1)+\frac{(2\mu+4)\left[3^{t}(\mu+1)^{t}-1\right]}{(\mu+1)(3\mu+2)}.
\end{equation}

By far, we finish derivations of Wiener index and mean first-passage time on three types of stochastic uniform growth trees. Compared to prior work focusing on deterministic versions, this study focuses on more general versions, and thus the formulas derived herein are also general. More importantly, the proposed method is more convenient to calculate what we care about than the commonly-used methods in the literature \cite{Agliari-2008,Zhang-2010,Wang-2012,Lin-2010}. In addition, we observe some differences and similarities between these growth trees using a systematical study, which is not yet reported in the previous work mainly because a single type of tree is often selected as an objective. This is helpful to understand the underlying structures on these growth trees. It is worth mentioning that during the derivation, a few surprising results are found, for instance, Eq.(\ref{Section-4-2-2-1}). This enables us to well reveal effect of graphic operations on topological structure of growth trees, and further create more intriguing networked models. Note also that while some other more complicated uniform growth trees are built and studied in the future, this work provides a guide to discuss many topological parameters including Wiener index and mean first-passage time on those models. Besides that, several example trees output by framework $\Upsilon$ can be selected to serve as candidate models modelling real-world networks \cite{Jurjiu-2003,Bartolo-2016,Furstenberg-2015,Markelov-2018}. Accordingly, the derived results are able to help one investigate topological structures on those networks.

\subsection{Network robustness and other structural parameters}

Network criticality $\overline{\mathcal{R}_{\mathcal{G}^{\dagger}}}$, as a topological measure estimating robustness on the underlying structure of network $\mathcal{G}(\mathcal{V},\mathcal{E})$ under consideration, has been widely studied in the past years \cite{Jekel-2018,Tizghadam-2010}. It is easy to see that in tree $\mathcal{T}$, quantity $\overline{\mathcal{R}_{\mathcal{T}^{\dagger}}}$ is in fact equal to average shortest path length $\overline{\mathcal{W}_{\mathcal{T}}}$. Therefore, we have the ability to analytically derive the corresponding solutions for such a parameter on all the stochastic uniform growth trees generated based on framework $\Upsilon$. Due to space limitation, we omit the correspondingly analytic expressions. On the other hand, we are interested in the scaling behavior of these parameters in the large graph size limit.

\textbf{Theorem 18} As $t\rightarrow\infty$, the asymptotic formulae for network criticality, $\overline{\mathcal{R}_{\mathcal{T}_{I}(t)^{\dagger}}}$, $\overline{\mathcal{R}_{\mathcal{T}_{II}(t)^{\dagger}}}$ and $\overline{\mathcal{R}_{\mathcal{T}_{III}(t)^{\dagger}}}$, on three classes of stochastic uniform growth trees are written as

\begin{equation}\label{Section-4-4-0}
\overline{\mathcal{R}_{\mathcal{T}_{I}(t)^{\dagger}}}=O(t),\quad \overline{\mathcal{R}_{\mathcal{T}_{II}(t)^{\dagger}}}=O\left(|V_{\mathcal{T}_{II}(t)}|^{\Theta_{II}-1}\right),\quad \text{and},\quad \overline{\mathcal{R}_{\mathcal{T}_{III}(t)^{\dagger}}}=O\left(|\mathcal{V}_{\mathcal{T}_{III}(t)}|^{\Theta_{III}-1}\right).
\end{equation}

These expressions are easily calculated and thus we omit proofs here. From the above equation, we can find that tree $\mathcal{T}_{I}(t)$ is more robust than other both types of stochastic uniform growth trees. One of most important reasons for this is that tree $\mathcal{T}_{I}(t)$ has a smaller diameter than trees $\mathcal{T}_{II}(t)$ and $\mathcal{T}_{III}(t)$.

In addition, trees $\mathcal{T}_{II}(t)$ and $\mathcal{T}_{III}(t)$ exhibit fractal structure, however such a phenomenon is not observed on tree $\mathcal{T}_{I}(t)$. As pointed in the literature \cite{Gennes-1982} (see Eq.(II.1) in \cite{Gennes-1982} for more details), for a fractal $\mathcal{G}(\mathcal{V},\mathcal{E})$ in question, there is an identity

$$\Theta_{\mathcal{G}}=\frac{2}{d_{\mathcal{G}}},$$
where $d_{\mathcal{G}}$ is the spectral dimension of graph $\mathcal{G}(\mathcal{V},\mathcal{E})$ and $\Theta_{\mathcal{G}}$ complies to $\overline{\mathcal{F}}_{\mathcal{G}}=|\mathcal{V}|^{\Theta_{\mathcal{G}}}$. So, using two parameters established above, $\Theta_{II}$ and $\Theta_{III}$, we can have

$$d_{\mathcal{T}_{II}(t)}=\frac{2\ln\left[1+\left(\sum_{i=1}^{\mu}p_{i}m_{i}\right)\left(\nu\sum_{j=1}^{\nu}q_{j}n_{j}+1\right)\right]}{\ln\left(\sum_{i=1}^{\mu}p_{i}m_{i}+1\right)+\ln\left[1+\left(\sum_{i=1}^{\mu}p_{i}m_{i}\right)\left(\nu\sum_{j=1}^{\nu}q_{j}n_{j}+1\right)\right]}<2,$$
and
$$d_{\mathcal{T}_{III}(t)}=\frac{2\ln\left(\mu\sum_{i=1}^{\mu}p_{i}m_{i}+1\right)}{\ln\left(\mu\sum_{i=1}^{\mu}p_{i}m_{i}+1\right)+\ln \left(2\sum_{i=1}^{\mu}p_{i}m_{i}+1\right)}<2.$$
In view of statement in \cite{Gennes-1982}, we point out that a walker originating from a designed vertex on tree $\mathcal{T}_{II}(t)$ will return back to the vertex almost surely over time because the corresponding spectral dimension is no more than $2$. The similar conclusion also holds for tree $\mathcal{T}_{III}(t)$.

\section{Conclusion}

To conclude, we consider random walks on tree networks and study some structural parameters of interest. First of all, we introduce three kinds of graphic operations, i.e., $VUGM$, $EUGM$ along with $MUGM$, and propose a principled framework $\Upsilon$. Based on this, we generate three families of stochastic uniform growth trees. As a consequence, some previously reported deterministic cases including $T$-graph and Vicsek fractal are contained into our framework completely. Next, we determine the analytic solution to mean first-passage time on stochastic uniform growth trees built. In view of an identity between Wiener index and mean first-passage time on tree given by Eq.(\ref{Section-3-4-0}), we first derive the corresponding formulae for Wiener index on stochastic uniform growth trees in a more manageable combinatorial manner instead of the commonly-used methods including spectral technique. One of most important reasons for this is that those typical manners mentioned above will become prohibitively difficult to execute even in some specific cases where the seed is just an edge or a star for instance. It should be emphasized that the formulae derived by us are more general, and thus cover the published results associated with deterministic cases. Last but not the least, we distinguish network robustness on all the stochastic uniform growth trees using network criticality. After that, we also analytically obtain spectral dimensions of two types of trees $\mathcal{T}_{II}(t)$ and $\mathcal{T}_{III}(t)$, and find a walker originating from a designed vertex on either $\mathcal{T}_{II}(t)$ or $\mathcal{T}_{III}(t)$ will return back to the vertex almost surely as time goes on.

We would like to stress that our work is only a tip of the iceberg and however the lights shed by our
methods can be beneficial to study random walks on other models \cite{Cohen-2016}-\cite{Alev-2020}. Meanwhile, there are still some open questions, for instance, how to effectively determine the $\mathcal{W}$-polynomial on the resulting graphs obtained from an arbitrary graph by utilizing three kinds of graphic operations introduced herein, which are waiting for us to address.

\section{Acknowledgments}

The authors would like to thank Xudong Luo for useful conversations. The research was supported by the National Key Research and Development Plan under grant 2020YFB1805400 and the National Natural Science Foundation of China under grant No. 61662066.

\section*{Appendix}

In Fig.3, we provide an illustrative example to clarify terminologies introduced in footnote 1.

\begin{figure}
\centering
  \includegraphics[height=4cm]{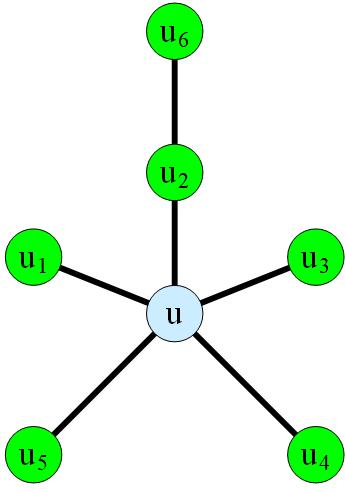}\\
  \vskip0.5cm
{\small Fig.3. (Color online) The diagram of a star-like graph centered at vertex $u$. More specifically, vertex $u$ is the center. There are five tentacles, namely, paths $\mathcal{P}_{uu_{1}}$, $\mathcal{P}_{uu_{3}}$, $\mathcal{P}_{uu_{4}}$, $\mathcal{P}_{uu_{5}}$ and $\mathcal{P}_{uu_{6}}$.}
\end{figure}

\end{document}